\documentclass[12pt,leqno]{article}

\usepackage[]{amsmath, amsthm, amsfonts}
\usepackage{amssymb}            
\usepackage[mathcal]{euscript}

\usepackage[matrix,arrow,curve,frame]{xy}       
\xymatrixcolsep{1.9pc}                          
\xymatrixrowsep{1.9pc}
\newdir{ >}{{}*!/-5pt/\dir{>}}

\numberwithin{equation}{section}

\newtheorem{theorem}{Theorem}[section]
\newtheorem{proposition}[theorem]{Proposition}
\newtheorem{lemma}[theorem]{Lemma}

\newtheorem{corollary}[theorem]{Corollary}
\newtheorem{definition}[theorem]{Definition}
\newtheorem{remark}[theorem]{Remark}

\newtheorem{notation}[theorem]{Notation}

\DeclareMathOperator{\Map}{Map}
\DeclareMathOperator{\Hom}{Hom}
\DeclareMathOperator{\Tot}{Tot}
\DeclareMathOperator{\TOT}{{\widetilde{Tot}}}
\DeclareMathOperator{\lan}{Lan}
\DeclareMathOperator{\hocolim}{hocolim}
\DeclareMathOperator{\colim}{colim}
\DeclareMathOperator{\N}{N}

\newcommand{\bx}{\mathbin{\square}}

\newcommand{\Proof}{\noindent{\bf Proof.\ }}
\newcommand{\id}{\mbox{\rm id}}
\renewcommand{\b}{\bullet}
\newcommand{\A}{{\cal A}}
\newcommand{\B}{{\cal B}}
\newcommand{\C}{{\cal C}}
\newcommand{\D}{{\cal D}}
\newcommand{\F}{{\cal F}}
\newcommand{\K}{{\cal K}}

\newcommand{\Q}{{\cal Q}}

\renewcommand{\O}{{\cal O}}
\newcommand{\QED}{\qed\bigskip}
\newcommand{\ov}{\overline}

\newlength{\labwidth}
\newcommand{\labarrow}[1]{
\settowidth{\labwidth}{$\scriptstyle \;\; #1 \;\;$}
\stackrel{#1}{\smash{\hbox to \labwidth{\rightarrowfill}} 
\vphantom{\longrightarrow}}
}

\begin{document}

\title{Cosimplicial objects and little $n$-cubes. I.}
\author{James E. McClure and Jeffrey H. Smith%
\thanks{Both authors were partially supported by NSF grants. 
The first author would also like to thank the
Isaac Newton Institute for its hospitality during the time this paper was being
written.}
\\Department of Mathematics, Purdue University  \\
150 N. University Street \\
West Lafayette, IN  47907-2067
}
\date{January 19, 2004}
\maketitle

\begin{abstract}
In this paper we show that if a cosimplicial space has a certain kind of 
combinatorial structure then its total space has an action of an operad weakly
equivalent to the little $n$-cubes operad.  Our results are also valid for 
cosimplicial spectra.
\end{abstract}

\section{Introduction.}

The little $n$-cubes operad $\C_n$ 
was introduced by Boardman and Vogt in
\cite{BV} (except that they used the terminology of theories rather than that 
of operads) as a tool for understanding $n$-fold loop spaces. 
They showed that for any topological space $Y$ the $n$-fold loop space 
$\Omega^n Y$ has an action of $\C_n$. 
In the other direction, May showed
in \cite{MayG} that if $Z$ is a space with an action of $\C_n$ then 
there exists a space $Y$ such that the group completion of $Z$ is weakly 
equivalent to $\Omega^n Y$.  

In the 30 years since \cite{MayG} the operad $\C_n$ has played an important 
role in both unstable and stable homotopy theory.  More recently, it has also 
been of importance (especially when $n=2$) in quantum algebra and other areas 
related to mathematical physics (see, for example, \cite{Get}, \cite{K1}, 
\cite{K2}, \cite{T1}, \cite{T2}).

In the known applications, $\C_n$ can be replaced by any operad weakly 
equivalent to it; such operads are called $E_n$ operads.

There is a highly developed technology that provides sufficient conditions for
a space to have an action by an $E_\infty$ operad (see \cite{Adams}, for
example) or an $E_1$ operad (\cite{MayG}, \cite{Thomason}).  Much less is known
about actions of $E_n$ operads for $1<n<\infty$. 

In this paper we consider the important special situation where the space (or 
spectrum) $Y$ on which we want an $E_n$ operad to act is obtained by 
totalization from a cosimplicial space (resp., spectrum) $X^\b$.  We 
construct an $E_n$ operad $\D_n$ and we show 
(Theorem \ref{6.4}) that if $X^\b$ has a certain kind of combinatorial 
structure (we call it a $\Xi^n$-structure) then $\D_n$ acts on $\Tot(X^\b)$.  

The converse of Theorem \ref{6.4} is not true: a $\D_n$-action on 
$\Tot(X^\b)$ does not have to come from a $\Xi^n$-structure on $X^\b$.  
However, in a future paper we will show that Tot induces a Quillen 
equivalence between the category of cosimplicial spaces with $\Xi^n$-structure 
and the category of spaces with $\D_n$-action; from this it will follow that 
if $\D_n$ acts on a space $Y$ then there is a cosimplicial space $X^\b$ with 
a $\Xi^n$ structure such that $\Tot(X^\b)$ is weakly equivalent to $Y$ as a 
$\D_n$-space.

The $n=2$ case of Theorem \ref{6.4} was proved in \cite{MS} (by a more 
complicated method than in the present paper), and it has had useful
applications, notably
\begin{itemize}
\item 
a proof (\cite{MS}) that the topological Hochschild cohomology spectrum 
of an $A_\infty$ ring spectrum $R$ has an action of $\D_2$ (this is the
topological analog of Deligne's Hochschild cohomology conjecture
\cite{Deligne}), and
\item
a proof by Dev Sinha that for $k\geq 4$ a space closely related to 
the space of knots in ${\mathbb R}^k$ 
is a 2-fold loop space. 
\end{itemize}

The methods we use in this paper are quite general and apply to other
categories of cosimplicial objects besides the categories of cosimplicial
spaces and spectra.
In the sequel to this paper \cite{MS4} we will apply our methods to the 
category of cosimplicial chain complexes. 

\begin{remark}
\label{F2}
{\rm
Throughout this paper we will use the following conventions for cosimplicial
spaces.

(a)
We define $\Delta$ to be the category of nonempty finite totally ordered
sets (this is equivalent to the category usually called $\Delta$).  We write 
$[m]$ for the finite totally ordered set $\{0,\ldots,m\}$.

(b) A cosimplicial space $X^\b$ is a functor from $\Delta$ to spaces.
If $S$ is a nonempty finite totally ordered set we write
$X^S$ for the value of $X^\b$ at $S$, except that we write $X^m$ instead of
$X^{[m]}$.
}
\end{remark}

Here is an outline of the paper.  

As an introduction to the ideas we begin in 
Sections \ref{sec2} and \ref{sec2a} with the $n=1$ (that is, the $A_\infty$) 
case.  In Section \ref{sec2} we recall the monoidal structure $\bx$ on the
category of cosimplicial spaces due to Batanin \cite{Bat93}.  In Section 
\ref{sec2a} we give a very simple proof of the fact (first shown in 
\cite{Bat98} and \cite{MS}) that if $X^\b$ is a monoid with respect to $\bx$ 
then $\Tot(X^\b)$ is an $A_\infty$ space; 
the proof is based on an idea due to Beilinson (\cite[Section 2]{HS}). 
We also give (in Remark \ref{W2}) an explicit 
description of the combinatorial structure on $X^\b$ that constitutes a 
$\bx$-monoid structure.

Our treatment of the $E_\infty$ case 
is precisely parallel, with
a symmetric
monoidal structure $\boxtimes$ in place of the monoidal structure $\bx$ (but
for technical reasons we need to use augmented cosimplicial spaces instead of
cosimplicial spaces; see Remark \ref{aug}).
As a prelude, in  Section \ref{sec3} we give a technically convenient
reformulation of the concept of symmetric monoidal structure; we also define
a more general concept (functor-operad)
which is used in later sections as a way of 
interpolating between monoidal and symmetric monoidal structures.
(The definition of functor-operad was discovered independently, in a different
context, by Batanin \cite{Bat02}).
In Section \ref{sec3a} we digress to offer motivation for the definition of
$\boxtimes$; the definition itself, and the verification that $\boxtimes$ is 
indeed a symmetric monoidal structure, is given in Section \ref{sec4}.
The main result in Section 
\ref{sec4a} (Theorem \ref{4.5}) is that
if $X^\b$ is a commutative monoid with respect to $\boxtimes$ then $\Tot(X^\b)$
is an $E_\infty$ space.  Section \ref{sec4a} also gives an explicit 
description of the combinatorial structure on $X^\b$ that constitutes a
commutative $\boxtimes$-monoid structure.

In Section \ref{sec5} we define for each $n$ a functor-operad $\Xi^n$; the
special case $n=\infty$ is the symmetric monoidal structure $\boxtimes$, and
the special case $n=1$ is essentially the same (see Remark \ref{MM2}) as the 
monoidal structure
$\bx$. (In the special case $n=2$, a construction isomorphic 
to $\Xi^2$ was discovered independently by Tamarkin in unpublished work.) 
In Section \ref{sec5a} we use the functor-operad $\Xi^n$ to construct
an ordinary (topological) operad $\D_n$.  The main theorem in Section
\ref{sec5a} (Theorem \ref{6.4}) says that $\D_n$ is weakly equivalent to 
the little $n$-cubes operad
$\C_n$ and that if $X^\b$
is a $\Xi^n$-algebra then $\D_n$ acts on $\Tot(X^\b)$.  Section \ref{sec5a} 
also gives an explicit description of the combinatorial structure on $X^\b$ that
constitutes an action of $\Xi^n$ on $X^\b$.

In Section \ref{sec6} we restrict to the special case $n=2$: we show that
an action of $\Xi^2$ is essentially the same thing as a more familiar
structure, namely an operad with multiplication as defined in \cite{GV}.
It seems likely that there is an 
analogous result for $n>2$, perhaps using Batanin's concept of $n$-operad 
\cite{Bat02}.

In Section \ref{sec7} we show that the operad $\D_n$ acts naturally on 
$\Omega^n Y$ for all $Y$.

The next two sections contain material which is used in the proofs of Theorems
\ref{4.5} and \ref{6.4} and may also be of independent interest.
Let $Y_k^0$ denote the $0$-th space of the cosimplicial space 
$(\Delta^\b)^{\boxtimes k}$.  In Section \ref{secA} we prove that 
$(\Delta^\b)^{\boxtimes k}$ is isomorphic as a cosimplicial space to
the Cartesian product $\Delta^\b\times Y_k^0$.  We also show that $Y_k^0$ has a
canonical cell structure and that it is contractible (which completes the
proof of Theorem \ref{4.5}). In
Section \ref{secB} we show that $\Xi^n_k(\Delta^\b,\ldots,\Delta^\b)$ is
isomorphic as a cosimplicial space to the Cartesian product of its $0$-th space
with $\Delta^\b$; this is used in the proof of Theorem \ref{6.4} (the $n=2$
case is the ``fiberwise prismatic subdivision'' used in \cite{MS}).

In Section \ref{secJ} we use a technique of Clemens Berger \cite{Berger} to 
show that the operad $\D_n$ defined in Section \ref{sec5a} is weakly 
equivalent to $\C_n$; this completes the proof 
of Theorem \ref{6.4}.  The basic idea of Berger's technique is to compare the
two operads by writing them as homotopy colimits, over the same indexing
category, of contractible operads.

Finally, in Section \ref{sec8} we observe that there is a variant of Tot
which preserves weak equivalences, and we show that Theorems \ref{4.5} and 
\ref{6.4} have analogs for this version of Tot.

We would like to thank the referee for a very careful reading of the paper 
and for several useful suggestions.

%\newpage

\section{A monoidal structure on the category of cosimplicial spaces.}

\label{sec2}

We begin with some motivation.  We are concerned with the question of 
when Tot of a cosimplicial space has an $A_\infty$ structure.  This question is
formally analogous to the question of when the normalization of a cosimplicial
abelian group has an $A_\infty$ structure (we will explore this analogy 
further in \cite{MS4}; also see \cite[Sections 3 and 4]{MSsurv}).
We take as our starting point the fact that for any space $W$, the normalized 
cochain complex of $W$ has a (strictly associative) product, namely the
cup product.  The normalized cochain complex of $W$ is the normalization of the
cosimplicial abelian group $S^\b W$ defined by 
\[
S^\b W=\Map (S_\b W, {\mathbb Z}),
\]
where $S_\b W$ is the singular complex of $W$ and Map
is set maps.   We therefore examine the relationship between the cup product
and the cosimplicial structure maps of $S^\b W$.

The cup product is defined as usual for $x\in S^p W$ and $y\in S^q W$
by
\[
(x\smallsmile y)(\sigma)=x(\sigma(0,\ldots,p))\cdot
y(\sigma(p,\ldots,p+q))
\]
Here $\sigma\in S_{p+q}W$, $\cdot$ is multiplication in $\mathbb Z$, and
$\sigma(0,\ldots,p)$ (resp., 
$\sigma(p,\ldots,p+q)$)
is the restriction of $\sigma$ to the subsimplex of
$\Delta^{p+q}$ spanned by the vertices $0,\ldots,p$ (resp., $p,\ldots,p+q$).
We note that $\smallsmile$ is related to the coface and codegeneracy 
operations by the following formulas:
\begin{equation}
\label{M1}
d^i (x\smallsmile y)=\left\{
\begin{array}{ll}
d^i x \smallsmile y & \mbox{if $i\leq p$} \\
x\smallsmile d^{i-p}y & \mbox{if $i>p$}
\end{array}
\right.
\end{equation}
\begin{equation}
\label{M2}
d^{p+1}x\smallsmile y=x\smallsmile d^0 y
\end{equation}
\begin{equation}
\label{M3}
s^i(x\smallsmile y)=\left\{
\begin{array}{ll}
s^i x\smallsmile y & \mbox{if $i\leq p-1$} \\
x \smallsmile s^{i-p}y & \mbox{if $i \geq p$}
\end{array}
\right.
\end{equation}

Now let us return to the category of cosimplicial spaces.  Formulas \eqref{M1},
\eqref{M2} and \eqref{M3} motivate the following definition.

\begin{definition}
\label{defbox}
{\rm 
Let $X^\b$ and $Y^\b$ be cosimplicial spaces. $X^\b\bx Y^\b$ is 
the cosimplicial space whose $m$-th space is 
\[
\left(\coprod_{p+q=m} X^p\times Y^q\right)/\sim
\]
(where $\sim$ is the equivalence relation generated by $(x,d^0 y)\sim 
(d^{|x|+1}x,y)$). 
The cosimplicial operators are given by
\[
d^i(x,y)=\left\{
\begin{array}{l}
(d^i x,y)\mbox{ if $i\leq |x|$} \\
(x,d^{i-|x|}y)\mbox{ if $i>|x|$}
\end{array}
\right.
\]
and 
\[
s^i(x,y)=\left\{
\begin{array}{l}
(s^i x,y)\mbox{ if $i\leq |x|-1$} \\
(x,s^{i-|x|}y)\mbox{ if $i\geq |x|$}
\end{array}
\right.
\]
}
\end{definition}

We leave it to the reader to check that the cosimplicial identities are
satisfied and that the following holds.

\begin{proposition}
$\bx$ is a monoidal structure for the category of
cosimplicial spaces, with unit the constant cosimplicial space that has a point
in every degree.  
\end{proposition}

There is also a monoidal structure $\bx$ for cosimplicial
spectra: one simply replaces the Cartesian products in Definition 
\ref{defbox} by smash products.

We conclude this section with some observations about Kan extensions.  This
material will not be needed logically for the rest of the paper, but it
provides useful motivation for the constructions in Sections 
\ref{sec4} and \ref{sec5}.

Recall the conventions in Remark \ref{F2}.

If $S_1,S_2,\ldots,S_k$ are finite totally ordered sets there is a unique 
total order 
on $S_1\coprod \cdots \coprod S_k$ for which the inclusion maps into the
coproduct are order-preserving and every element of $S_i$ is less than every 
element of $S_j$ for $i<j$.
Let
\[
\Phi: \Delta^{\times k}\rightarrow\Delta
\]
be the functor which takes $(S_1,\ldots,S_k)$ to 
$S_1\coprod \cdots \coprod S_k$ with this total order.

The following fact was first noticed by Cordier and Porter (unpublished). 

\begin{proposition}
\label{Kan}
Let $X_1^\b,\ldots,X_k^\b$ be cosimplicial spaces and let 
$X_1^\b\,\bar{\times}
\cdots\bar\times\,X_k^\b$ denote the composite
\[
\Delta^{\times k} \labarrow{X_1^\b\times \ldots\times X_k^\b}
\mbox{\rm Top}\times\cdots\times\mbox{\rm Top} \labarrow{\times} \mbox{\rm 
Top}.
\]
Then
$X_1^\b\bx \cdots \bx X_k^\b$ is naturally isomorphic to the left 
Kan extension 
$\mbox{\rm Lan}_\Phi (X_1^\b\,\bar\times\cdots\bar\times \,X_k^\b)$.
\end{proposition}

Before giving the proof we mention an important consequence.  

\begin{remark}
\label{W1}
{\rm
Let $\Phi^*$ be
the functor from cosimplicial spaces to bicosimplicial spaces defined by 
\[
(\Phi^*(X^\b))^{S,T}=X^{S\coprod T}
\]
It is a general fact about Kan extensions 
\cite[beginning of Section X.3]{MacLane}
that 
$\mbox{\rm Lan}_\Phi$ is the left adjoint of $\Phi^*$.  This implies that 
there is a natural 1-1 
correspondence between maps
\[
\alpha: X^\b \bx Y^\b \to Z^\b
\]
and consistent collections of maps
\[
\widetilde\alpha_{S,T}:
X^S \times Y^T \to Z^{S \coprod T}
\]
where ``consistent'' means that every pair of ordered maps $f:S\to S'$, 
$g:T\to T'$ induces a commutative diagram
\[
\xymatrix{
X^S \times Y^T \
\ar[d]_{f_*\times g_*}
\ar[r]^-{\widetilde\alpha_{S,T}}
&
\ Z^{S\coprod T}
\ar[d]^{(f\coprod g)_*}
\\
X^{S'} \times Y^{T'} \
\ar[r]^-{\widetilde\alpha_{S',T'}}
&
\ Z^{S'\coprod T'}
}
\]
}
\end{remark}

\noindent
{\bf Proof of Proposition \ref{Kan}.\ }
For simplicity we assume $k=2$.

For $m\geq 0$ let $[m]$ denote the set $\{0,\ldots,m\}$.
Every object in $\Delta$ is canonically isomorphic to one of the form
$[m]$ so it suffices to show that the two functors in question are naturally 
isomorphic on the full subcategory of $\Delta$ with these objects.

Fix $m$.
According to \cite[Equation (10) on page 240]{MacLane}, the left Kan extension,
evaluated at $[m]$, can be calculated as follows. Let $\cal C$ be 
the category of objects $\Phi$-over $[m]$: an object in $\cal C$ is a 
pair consisting of an object $([n],[n'])$ in $\Delta\times\Delta$ and a
morphism $\Phi([n],[n'])\rightarrow [m]$ in $\Delta$; a morphism in $\cal
C$ is a morphism $(f,g)$ in $\Delta\times \Delta$ making the evident triangle
commute.  
Then
$\mbox{Lan}_\Phi (X_1^\b\bar\times X_2^\b)$ is the colimit of the composite
\[
{\cal C}\rightarrow \Delta\times\Delta \labarrow{X_1^\b\bar\times X_2^\b} 
\mbox{ \rm Top}
\]
where the first map is the evident forgetful functor.
Now let ${\cal C}'$ be the full subcategory of $\cal C$ consisting of pairs
$(([p],[q]),f:[p+q+1]\rightarrow [m])$ where $f$ is a surjection, $p+q$ is 
either $m-1$ or $m$, and if $p+q=m$ then $f(p)=f(p+1)$.  Note that an object 
in ${\cal C}'$ is
determined by $p$ and $q$; to simplify the notation we denote the object by 
$(p,q)$.  If $p+q=m-1$ there is exactly one morphism from the object $(p,q)$ to
$(p+1,q)$ and exactly one from $(p,q)$ to $(p,q+1)$, and there are no other
non-identity morphisms in ${\cal C}'$.  The map from $(p,q)$ to $(p+1,q)$ is
the last coface map on $[p]$ and the identity on $[q]$, while the map from
$(p,q)$ to $(p,q+1)$ is the identity on $[p]$ and the zeroth coface map on
$[q]$.  From this it is clear that the $m$-th space of
$X_1^\b\bx X_2^\b$ is the colimit of the composite
\[
{\cal C}'\subset{\cal C}\rightarrow \Delta\times\Delta 
\labarrow{X_1^\b\bar\times X_2^\b} \mbox{ \rm Top}
\]
In particular there is a natural map 
\[
X_1^\b\bx X_2^\b\rightarrow 
\mbox{Lan}_\Phi (X_1^\b\bar\times X_2^\b).
\]  By \cite[Section 9.3]{MacLane}, this 
will be an
isomorphism if ${\cal C}'$ is cofinal in ${\cal C}$, that is, if for each $c\in
{\cal C}$ the under-category $c\downarrow {\cal C}'$ is connected.  The
under-categories can be described explicitly: each is a nonempty full 
subcategory of ${\cal C}'$ with set of objects of the form 
\[
\{(p,q) \,|\, p\geq p_0, q\geq q_0\}
\]
for some $p_0$ and $q_0$;
clearly all such categories are connected.

\QED

%\newpage

\section{A sufficient condition for $\Tot(X^\b)$ to be an $A_\infty$ space.}

\label{sec2a}

In this section we will prove

\begin{theorem} \label{A-infty}
If $X^\b$ is a monoid with respect to $\bx$ then $\Tot(X^\b)$ has an
$A_\infty$ structure.
\end{theorem}

\begin{remark} 
\label{T1}
{\rm 
(a) Previous proofs of Theorem \ref{A-infty} were given by Batanin 
\cite[Theorems 5.1 and 5.2]{Bat98} and by us \cite[Theorem 2.4]{MS}.

(b) The theorem and its proof are also valid for cosimplicial spectra.
}
\end{remark}

\begin{remark}
\label{W2}
{\rm
Definition \ref{defbox} implies that $X^\b$ is a monoid with respect to
$\bx$ if and only if there are maps
\[
\smallsmile: X^p\times X^q \to X^{p+q}
\]
for all $p,q\geq 0$ satisfying equations \eqref{M1}, \eqref{M2}, \eqref{M3},
the associativity condition
\begin{equation}
\label{F1}
(x\smallsmile y)\smallsmile z=z\smallsmile (y\smallsmile z)
\end{equation}
and the unit condition: there is an element $e\in X^0$ such that
\begin{equation}
\label{T2}
x\smallsmile e=e\smallsmile x=x
\end{equation}
for all $x$.
}
\end{remark}

The rest of this section is devoted to the proof of Theorem \ref{A-infty}.

Let $\Delta^\b$ denote the cosimplicial space whose $m$-th space is
the simplex $\Delta^m$, with the usual cofaces and codegeneracies.  By 
definition, $\Tot(X^\b)$ is $\Hom(\Delta^\b,X^\b)$ (where Hom denotes the space
of cosimplicial maps).

\begin{definition}
\label{T3}
{\rm
(a)
For each $k\geq 0$ let 
$\A(k)$ be the space $\Tot((\Delta^\b)^{\bx k})$.

(b)
If $f\in \A(k)$ and $g_i\in \A(j_i)$ for $1\leq i \leq k$ define
$\gamma(f,g_1,\ldots,g_k)\in \A(j_1+\cdots+j_k)$ 
to be the composite
\[
\Delta^\b\labarrow{f}(\Delta^\b)^{\bx k}
\labarrow{g_1\bx\cdots\bx g_k} 
(\Delta^\b)^{\bx(j_1+\cdots+j_k)}
\]
}
\end{definition}

Theorem \ref{A-infty} is an immediate consequence of our next result.  

\begin{proposition}
\label{2.3}
(a) Let $\A$ be the sequence of spaces $\A(k)$, $k\geq 0$, with
the operations 
\[
\gamma:\A(k)\times\A(j_1)\times\cdots\times\A(j_k)\rightarrow
\A(j_1+\cdots+j_k)
\]
defined above.
Then $\A$ is an operad.

(b) If $X^\b$ is a monoid with respect to $\bx$ then $\A$ acts
on $\Tot(X^\b)$.

(c) $\A$ is an $A_\infty$ operad.
\end{proposition}

\Proof  Part (a) is clear. 

For part (b), given $f\in
\A(k)$ and $x_1,\ldots,x_k\in\Tot(X^\b)$ define
$f(x_1,\ldots,x_k)\in\Tot(X^\b)$ to be the composite
\[
\Delta^\b\labarrow{f}(\Delta^\b)^{\bx k} 
\labarrow{x_1\bx\cdots\bx x_k} (X^\b)^{\bx k}
\labarrow{\mu} X^\b
\]
where $\mu$ is the monoidal structure map of $X^\b$.  This construction gives
maps
\[
\A(k)\times (\mbox{Tot}(X^\b))^k\rightarrow \mbox{Tot}(X^\b)
\]
which fit together to give an action of $\A$ on Tot$(X^\b)$.

For part (c) we need to show that $\A(0)$ is a point (which is obvious) and 
that each space $\A(k)$ is contractible.  First consider the case $k=1$.  If
$f$ and $g$ are two cosimplicial maps from $\Delta^\b$ to 
$\Delta^\b$, then
$tf+(1-t)g$ will again be a cosimplicial map for each $0\leq t \leq 1$ (because
the cosimplicial structure maps of $\Delta^\b$ are affine) and so we can use
the straight-line homotopy to contract Hom$(\Delta^\b,\Delta^\b)$ to a point.
The case $k\geq 2$ is now immediate from Lemma \ref{iso} below.
\QED

\begin{lemma} \label{iso}
For each $k\geq 1$, $(\Delta^\b)^{\bx k}$ is isomorphic as a 
cosimplicial
space to $\Delta^\b$.
\end{lemma}

\Proof It suffices to do the case $k=2$; the general case follows by
induction.  First we define maps 
\[
f^m: (\Delta^\b\bx\Delta^\b)^m\rightarrow \Delta^m
\]
for $m\geq 0$ by
\[
\textstyle
f^m((s_0,\ldots,s_p),(t_0,\ldots,t_q))=
\left(\frac{1}{2}s_0,\ldots,\frac{1}{2}(s_p+t_0),\ldots,\frac{1}{2}t_q\right)
\]
These are well defined and fit together to give a
cosimplicial map $f:\Delta^\b\bx\Delta^\b\rightarrow\Delta^\b$.  Next 
define
\[
g^m:\Delta^m\rightarrow(\Delta^\b\bx\Delta^\b)^m
\]
as follows: given $(u_0,\ldots,u_m)\in \Delta^m$, choose the smallest $p$ for
which 
\[
u_0+\cdots+u_p\geq \frac{1}{2}
\]
and let
\[
g^m(u_0,\ldots,u_m)=[(2u_0,\ldots,1-2u_0-\cdots-2u_{p-1}),
(2u_0+\cdots+2u_p-1,2u_{p+1},\ldots,2u_m)]
\]
The $g^m$ are continuous, they fit together to 
give a cosimplicial map $g:\Delta^\b\rightarrow \Delta^\b\bx\Delta^\b$,
and $f$ and $g$ are mutually inverse.
\QED

\begin{remark}
{\rm
Lemma \ref{iso} is due to Grayson \cite[Section 4]{G}.  
}
\end{remark}

%\newpage

\section{Functor-operads.}

\label{sec3}

The purpose of this section is to describe a general setting in which there are
analogs of Definition \ref{T3} and Proposition \ref{2.3}(a) and (b) (see
Definition \ref{extra} and Propositions \ref{3.2} and \ref{3.4}).

Given a category $\C$ let
$\C^{\times k}$ denote the $k$-fold Cartesian product. For
each permutation $\sigma\in \Sigma_k$ we define
\[
\sigma_\#: \C^{\times k} \rightarrow \C^{\times k}
\] 
to be the functor taking
$(A_1,\ldots,A_k)$ to $(A_{\sigma(1)},\ldots,A_{\sigma(k)})$.

In order to motivate the definition of functor-operad, let us consider 
the situation in which 
$\C$ has a symmetric monoidal structure $\boxtimes$.  For each
$k\geq 0$ define a functor 
\[
\F_k:\C^{\times k} \to \C
\]
by
\[
\F_k(X_1,\ldots,X_k)=X_1\boxtimes(X_2\boxtimes(X_3\boxtimes\cdots))
\]
MacLane's coherence theorem \cite{MacLane} implies that there are canonical
natural isomorphisms
\[
\sigma_*:\F_k\rightarrow \F_k\circ\sigma_\#
\]
and
\[
\Gamma_{j_1,\ldots,j_k}:
\F_k(\F_{j_1},\ldots,\F_{j_k})\rightarrow \F_{j_1+\cdots j_k}
\]
satisfying certain consistency conditions.  The following definition is an
abstract version of this situation, except that instead of requiring the
$\Gamma$'s to be natural isomorphisms we allow them merely to be natural
transformations.

\begin{definition}
\label{funop}
{\rm
Let $\C$ be a category enriched over Top.  A {\it functor-operad} $\F$ in 
$\C$ is a sequence of continuous
functors $\F_k:\C^{\times k}\rightarrow \C$ 
together with 

(i) for each $\sigma\in\Sigma_k$, a continuous natural isomorphism
\[
\sigma_*:\F_k\rightarrow \F_k\circ\sigma_\#
\]

(ii) for each choice of $j_1,\ldots,j_k\geq 0$, a continuous natural 
transformation
\[
\Gamma_{j_1,\ldots,j_k}: 
\F_k(\F_{j_1},\ldots,\F_{j_k})\rightarrow \F_{j_1+\cdots j_k}
\]

such that

(a) $\F_1$ is the identity functor, and the natural transformations
\[
\Gamma_{1,\ldots,1}:\F_k(\F_1,\ldots,\F_1)\rightarrow \F_k
\]
\[
\Gamma_{k}:\F_1(\F_k)\rightarrow \F_k
\]
are equal to the identity.

(b) All diagrams of the following form commute:
\[
\xymatrix{
\F_k(\F_{j_1}(\F_{i_{11}},\ldots,\F_{i_{1j_1}}),\ldots,
\F_{j_k}(\F_{i_{k1}},\ldots,\F_{i_{kj_k}})) 
\ar[r]^-{\Gamma}
\ar[d]_{\F_k(\Gamma,\ldots,\Gamma)}
&
\F_{j_1+\ldots+j_k}(\F_{i_{11}},\ldots,\F_{i_{kj_k}}) 
\ar[d]_{\Gamma}
\\
\F_k(\F_{i_{11}+\cdots+i_{1j_1}},\ldots,\F_{i_{k1}+\cdots+i_{kj_k}}) 
\ar[r]^-{\Gamma}
&
\F_{i_{11}+\cdots+i_{kj_k}}
}
\]

(c) $(\sigma\tau)_*=\tau_*\sigma_*$ for all $\sigma, \tau\in \Sigma_k$

(d) Let $\tau_i\in\Sigma_{j_i}$ for $1\leq i\leq k$ and let
$\tau$ be the image of $(\tau_1,\ldots,\tau_k)$ under the map 
\[
\Sigma_{j_1}\times\cdots\times\Sigma_{j_k}\rightarrow \Sigma_{j_1+\cdots+j_k}
\]
Then the following diagram commutes:
\[
\xymatrix{
\F_k(\F_{j_1},\ldots,\F_{j_k}) 
\ar[d]_{\F_k(\tau_{1*},\ldots,\tau_{k*})}
\ar[r]^-{\Gamma}
&
\F_{j_1+\cdots+j_k} 
\ar[d]_{\Gamma}
\\
\F_k(\F_{j_1}\circ \tau_{1\#},\ldots,\F_{j_k}\circ \tau_{k\#}) 
\ar[r]^-{\tau_*}
& 
\F_{j_1+\cdots+j_k}\circ \tau_\#
}
\]

(e) Let $\sigma\in\Sigma_k$ and let  $\bar{\sigma}$ be the permutation in
$\Sigma_{j_1+\ldots+j_k}$ which permutes the blocks
$\{1,\ldots,j_1\},\ldots,\{j_1+\ldots+j_{k-1}+1,\ldots,j_1+\ldots+j_k\}$
in the same way that $\sigma$ permutes the numbers $1,\ldots,k$.
Then the following diagram commutes
\[
\xymatrix{
\F_n(\F_{j_1},\ldots,\F_{j_k}) 
\ar[r]^-{\Gamma}
\ar[d]_{\sigma_*}
&
\F_{j_1+\cdots+j_k} 
\ar[d]_{\bar{\sigma}_*}
\\
\F_n(\F_{j_{\sigma(1)}},\ldots,\F_{j_{\sigma(k)}})\circ \bar{\sigma}_\# 
\ar[r]^-{\Gamma}
&
\F_{j_1+\cdots+j_k}\circ \bar{\sigma}_\#
}
\]
}
\end{definition}

\begin{remark}
\label{j}
{\rm
(a) Batanin \cite{Bat02} has independently proposed a similar but more general 
definition: if $\O$ is an operad in the category of categories Batanin 
defines an {\it internal operad} in $\O$ to be a collection consisting of an 
object $a_k$ in $\O(k)$ for each $k\geq 0$ and morphisms 
\[
\sigma_*: a_k\rightarrow \tilde{\sigma}(a_k)
\]
for each $\sigma\in\Sigma_k$
(where $\tilde{\sigma}$ denotes the action of $\sigma\in\Sigma_k$ on 
$\O(k)$) and
\[
\Gamma_{j_1,\ldots,j_k}:
\gamma(a_k,a_{j_1},\ldots,a_{j_k})\rightarrow a_{{j_1}+\cdots+a_{j_k}}
\]
for each $j_1,\ldots,j_k\geq 0$ (where $\gamma$ is the structure map of the
operad $\O$) satisfying the analogs of properties (a)--(e) in Definition 
\ref{funop}.  A functor-operad in $\C$ is then an internal operad in the
endomorphism operad of $\C$.

(b) If $\B$ is an operad in the category Top we can define a functor-operad
$\F$ in Top by 
\[
\F_k(X_1,\ldots,X_k)=\B_k\times X_1\times\cdots\times X_k
\]
with the obvious structure maps.
}
\end{remark}

\medskip

\begin{definition}
{\rm
\label{extra}
Let $\F$ be a functor-operad in $\C$ and let $A$ be an object of $\C$. 

(a)
Define $\F_A$ to be the collection of spaces 
\[
\F_A(k)=\Hom(A,\F_k(A,\ldots,A)), \quad k\geq
0.
\] 

(b)
Give $\F_A(k)$ the action induced by the $\sigma_*$.

(c)
Define $1\in \F_A(1)$ to be the identity map of $A$.

(d)
For each
choice of $j_1,\ldots,j_k\geq 0$ define
\[
\gamma: 
\F_A(k) \times
\F_A(j_1)\times  
\cdots
\F_A(j_k)
\rightarrow
\F_A(j_1+\cdots+j_k)
\]
to be the composite 
\[
\begin{array}{l}
\Hom(A,\F_k(A,\ldots,A))\times  
\Hom(A,\F_{j_1}(A,\ldots,A))\times  
\cdots
\Hom(A,\F_{j_k}(A,\ldots,A))
\rightarrow \\
\quad\Hom(A,\F_k(\F_{j_1}(A,\ldots,A),\ldots,\F_{j_k}(A,\ldots,A)))
\labarrow{{\scriptstyle\mathrm H\mathrm o\mathrm m}(A,\Gamma)}
\Hom(A,\F_{j_1+\cdots+j_k}(A,\ldots,A))
\end{array}
\]
}
\end{definition}

\begin{proposition}
\label{3.2}
These choices make $\F_A$ an operad.
\end{proposition}

The proof is an easy verification.

\begin{definition}
\rm{
Let $\F$ be a functor-operad in $\C$.  An {\it algebra over $\F$} is an object
$X$ of $\C$ together with continuous maps 
$\Theta_k:\F_k(X,\ldots,X)\rightarrow X$ for $k\geq 0$ such that

(a) $\Theta_1$ is the identity map. 

(b) The following diagram commutes for each choice of $j_1,\ldots,j_k\geq 0$
\[
\xymatrix{
\F_k(\F_{j_1}(X,\ldots,X),\ldots,\F_{j_k}(X,\ldots,X)) 
\ar[r]^-{\Gamma}
\ar[d]_{\F_k(\Theta_{j_1},\ldots,\Theta_{j_k})}
&
\F_{j_1+\cdots+j_k}(X,\ldots,X)
\ar[d]^{\Theta_{j_1+\cdots+j_k}}
\\
\F_k(X,\ldots,X)
\ar[r]^-{\Theta_k}
&
X
}
\]

(c) $\Theta_k\circ \sigma_*=\Theta_k$ for all $\sigma\in \Sigma_k$. 

}
\end{definition}

Now let $A$ be an object of $\C$, let $X$ be an algebra over $\F$, and for each
$k\geq 0$ define
\[
\theta_k:\F_A(k)\times \Hom(A,X)^k\rightarrow \Hom(A,X)
\]
to be the composite
\[
\Hom(A,\F_k(A,\ldots,A))\times \Hom(A,X)^k \rightarrow 
\Hom(A,\F_k(X,\ldots,X)) \labarrow{\mathrm H \mathrm o \mathrm m (A,\Theta_k)} 
\Hom(A,X)
\]

\begin{proposition}
\label{3.4}
The maps $\theta_k$ make $\Hom(A,X)$ an algebra over the operad $\F_A$.
\end{proposition}

Again, the proof is an easy verification.

\begin{definition}
{\rm
A functor-operad $\F$ is {\it strict} if the natural transformations
$\Gamma_{j_1,\ldots,j_k}$ are isomorphisms.
}
\end{definition}

\begin{proposition}
\label{3.6}
If $\F$ is a strict functor-operad, then $\F_2$ is a symmetric monoidal 
structure for $\C$ with identity object $\F_0$.  The commutative monoids with
respect to this structure are the same as the algebras over $\F$.
\end{proposition}

Once more, the proof is an easy verification.

%\newpage

\section{A family of operations in $S^\b W$}

\label{sec3a}

In Section \ref{sec4} we will define a symmetric monoidal product $\boxtimes$
on the category of augmented cosimplicial spaces. 
In this section we pause to offer motivation for this definition.  The results
in this section are not needed logically for later sections.

The definition of the monoidal product $\bx$ 
was motivated in Section \ref{sec2} by the
properties of the cup product in $S^\b W$.  
The cup product is part of a larger family of operations in $S^\b W$ whose 
properties could be used as the basis for a definition of $\boxtimes$.
However, this larger family is rather inconvenient to work with
(because the analog of equation \eqref{M2} for the larger family is
complicated) so we will use a related family which has somewhat
simpler properties. 

We begin with a variant of the cup 
product.  Given $x\in S^p W$ and $y\in S^q W$
we define
\[
x\sqcup y\in S^{p+q+1} W
\]
by
\[
(x\sqcup y)(\sigma)=x(\sigma(0,\ldots,p))\cdot
y(\sigma(p+1,\ldots,p+q))
\]
(note that, in contrast to the cup product, the vertex $p$ is not repeated). 

This operation is related to the coface and codegeneracy operations in $S^\b W$
by the following equations:
\begin{equation}
\label{T4}
d^i (x\sqcup y) =
\left\{
\begin{array}{ll}
d^i x \sqcup y & \mbox{if $i\leq p+1$} \\
x\sqcup d^{i-p-2}y & \mbox{if $i>p+1$}
\end{array}
\right.
\end{equation}
\begin{equation}
\label{T5}
s^i (x\sqcup y) =
\left\{
\begin{array}{ll}
s^i x \sqcup y & \mbox{if $i< p$} \\
x\sqcup s^{i-p-1}y & \mbox{if $i>p$}
\end{array}
\right.
\end{equation}
Note that there is no analog for $\sqcup$ of equation \eqref{M2}.

The operations $\smallsmile$ and $\sqcup$ determine each other:
\[
x\sqcup y = (d^{p+1}x)\smallsmile y =x\smallsmile d^0 y
\]
\[
x \smallsmile y =s^p(x \sqcup y)
\]

Now observe that equations \eqref{T4} and \eqref{T5} can be used as the
basis for a characterization of $\bx$-monoids: Remark \ref{W1} implies that 
$X^\b$ is a 
$\bx$-monoid if and only if there are maps
\[
\sqcup:X^p\times X^q\to X^{p+q+1}
\]
satisfying \eqref{T4}, \eqref{T5}, the associativity condition
\begin{equation}
\label{T6}
x\sqcup (y\sqcup z)=
(x\sqcup y)\sqcup z
\end{equation}
and the unit condition: there exists $e\in X^0$ with
\begin{equation}
\label{T7}
s^p(x\sqcup e)=s^0(e\sqcup x)=x
\end{equation}
(compare this to Remark \ref{W2}).

In the remainder of this section we will define a family of operations in
$S^\b W$ which generalize $\sqcup$; the definition of $\boxtimes$ in Section 
\ref{sec4} will be suggested by the properties of this family.

First we need to be a little more explicit in our description of $S^\b W$.  
Recall the conventions in Remark \ref{F2}.  Given a nonempty finite set totally
ordered set $T$ we let $\Delta^T$ be the convex hull of $T$; in
particular, $\Delta^{[m]}$ is the usual $\Delta^m$.  We define $S_T W$ to be
the set of all continuous maps $\Delta^T\to W$ (in particular, $S_{[m]} W$ is
what we have been calling $S_m W$) and $S^T W$ to be
$\Map(S_T W, {\mathbb Z})$ (so $S^{[m]} W$ is the same as $S^m W$).

Given a map
$\sigma:\Delta^T\to W$ and a subset $U$ of $T$ let $\sigma(U)$
denote the restriction of $\sigma$ to the sub-simplex of $\Delta^T$ spanned by
the vertices in $U$.  

Suppose we are given a function
\[
f:T \to \{1,\ldots,k\}
\]
and elements $x_i\in S^{f^{-1}(i)}W$ for $1\leq i\leq k$. We can define 
an element
\[
\langle f \rangle(x_1,\ldots,x_k) \in S^T W
\]
by
\[
\langle f \rangle(x_1,\ldots,x_k)(\sigma)=
x_1(\sigma(f^{-1}(1))\cdot
x_2(\sigma(f^{-1}(2)) \cdot
\,\cdots\,
\cdot x_k(\sigma(f^{-1}(k))
\]
where $\cdot$ denotes multiplication in $\mathbb Z$.  This procedure gives a
natural transformation
\[
\langle f \rangle:
S^{f^{-1}(1)}W \otimes
\cdots \otimes
S^{f^{-1}(k)}W
\to
S^T W
\]

\begin{remark}
\label{F4}
{\rm
In the special case where
$f$ is the function from $\{0,\ldots,p+q+1\}$ to $\{1,2\}$ which takes
$\{0,\ldots,p\}$ to 1 and $\{p+1,\ldots,p+q+1\}$ to 2, we have 
$
\langle f\rangle (x,y)= x\sqcup y
$
}
\end{remark}

Next we describe the relation between the operations $\langle f\rangle$ and 
the cosimplicial structure maps of $S^\b W$. 

\begin{proposition} 
\label{F3}
Let 
\[
\xymatrix{
T
\ar[rr]^h
\ar[dr]_-f
&&
T'
\ar[dl]^-g
\\
&
\{1\ldots,k\}
&
}
\]
be a commutative diagram, where 
$h$ is a map in $\Delta$ (i.e., an order-preserving map).  For each $i\in
\{1\ldots,k\}$ let 
\[
h_i:f^{-1}(i)\to g^{-1}(i)
\]
be the restriction of $h$.  

Then the diagram
\[
\xymatrix{
S^{f^{-1}(1)}W\otimes\cdots\otimes S^{f^{-1}(k)}W
\ar[d]_{(h_1)_*\otimes\cdots\otimes (h_k)_*}
\ar[r]^-{\langle f\rangle}
&
S^T W
\ar[d]^{h_*}
\\
S^{g^{-1}(1)}W\otimes\cdots\otimes S^{g^{-1}(k)}W
\ar[r]^-{\langle g\rangle}
&
S^{T'}W
}
\]
commutes.  
\end{proposition}

The proof is an immediate consequence of the definitions.  In the special 
case of Remark \ref{F4} we recover equations \eqref{T4} and \eqref{T5}.

%\newpage

\section{A symmetric monoidal structure on the category of augmented 
cosimplicial spaces.}

\label{sec4}

From now on we will work with augmented cosimplicial spaces (the reason for 
this is given in Remark \ref{aug}).  

\begin{definition}
{\rm
An augmented cosimplicial
space is a functor $X^\b$ from $\Delta_+$ to Top, where $\Delta_+$ is the 
category of finite totally ordered sets (including the empty set).
}
\end{definition}

Our goal in this section is to construct a symmetric monoidal product
$\boxtimes$ in the category of augmented cosimplicial spaces.  We will do 
this by constructing a strict functor-operad $\Xi$ and letting $\boxtimes$ be 
$\Xi_2$; see Proposition \ref{3.6}. 

The basic idea in defining
$\Xi_k(X_1^\b,\ldots,X_k^\b)$ is that we build it from formal symbols
$\langle f \rangle (x_1,\ldots,x_k)$, where $f:T\to \{1,\ldots,k\}$ (cf.\ 
Section \ref{sec3a}).  In order 
to get a cosimplicial object we have to
build in the cosimplicial operators, so we consider symbols of the form
\[
h_*(\langle f \rangle (x_1,\ldots,x_k))
\]
where $h:T\to S$ is an order-preserving map; such a
symbol will represent a point in the $S$-th space 
$\Xi_k(X_1^\b,\ldots,X_k^\b)^S$.
We want to require these symbols to satisfy the relation in Proposition 
\ref{F3}, and the most efficient way to do this is by means of a Kan extension
(Definition \ref{Kan2}); a more elementary description of $\Xi_k$ is given in 
equation \ref{F6}.

Here are the formal definitions:

\begin{definition}
\label{bark}
{\rm
Let $k\geq 0$. Define $\bar{k}$ to be
the set $\{1,\ldots,k\}$ when $k\geq 1$ and the empty set when $k=0$.
}
\end{definition}

\begin{definition}
\label{Q}
{\rm
Let $\Q_k$ be the category 
whose objects are pairs $(f,S)$, where $S$ is an object of $\Delta_+$ and $f$ 
is a map of sets from $S$ to $\bar{k}$, and whose morphisms 
are commutative triangles
\[
\xymatrix{
S
\ar[rr]^h
\ar[dr]_-f
&&
T
\ar[dl]^-g
\\
&\bar{k}&
}
\]
where $h$ is a map in $\Delta_+$.
}
\end{definition}

There is a forgetful functor $\Phi:\Q_k\rightarrow \Delta_+$ which takes 
$(f,S)$ 
to $S$, and a functor $\Psi$ from $\Q_k$ to the $k$-fold Cartesian product 
$(\Delta_+)^{\times k}$ which takes
$(f,S)$ to the $k$-tuple $(f^{-1}(1),\ldots,f^{-1}(k))$.

\begin{definition}
\label{Kan2}
{\rm
For each $k\geq 0$ define a functor $\Xi_k$ as follows.
Given augmented cosimplicial spaces $X_1^\b,\ldots,X_k^\b$, 
let $X_1^\b\bar\times\cdots\bar\times X_k^\b$ denote the composite
\[
(\Delta_+)^{\times k} \labarrow{X_1^\b\times\cdots\times X_k^\b}
\mbox{\rm Top} \times
\cdots\times
\mbox{\rm Top} 
\labarrow{\times} \mbox{\rm Top}.
\]
We define the augmented cosimplicial space $\Xi_k(X_1^\b,\ldots,X_k^\b)$ to 
be the Kan extension 
\[
\mbox{Lan}_\Phi ((X_1^\b\bar\times\cdots\bar\times X_k^\b)\circ \Psi)
\]
%\[
%\xymatrix{
%\Q(k) 
%\ar[d]_{\Phi}
%\ar[r]^{\Psi}
%& (\Delta_+)^{\times k} 
%\ar[rr]^{X_1^\b\bar\times\cdots\bar\times X_k^\b}
%&& \mbox{Top} \\
%\Delta_+ 
%\ar[rrru]|{\Xi_k(X_1^\b,\ldots,X_k^\b)}
%& &  
%}
%\]
}
\end{definition}

\begin{remark}
\label{r2}
{\rm
(a) $\Xi_0$ is the augmented cosimplicial space which takes every $S$ to a 
point (because a Cartesian product indexed by the empty set is a point).

(b) 
The adjointness property of $\lan_\Phi$
\cite[beginning of Section X.3]{MacLane} implies that a map 
$\Xi_k(X_1^\b,\ldots,X_k^\b)\to Y^\b$ is the same thing as a collection of maps
\[
\langle f \rangle: X_1^{f^{-1}(1)} \times \cdots \times X_k^{f^{-1}(k)}
\to Y^T,
\]
one for each $f:T\to \bar{k}$, such that the analog of Proposition
\ref{F3} is satisfied.
}
\end{remark}

Our next goal is to specify the structure maps $\sigma_*$ and 
$\Gamma_{j_1,\ldots,j_k}$ of the functor-operad $\Xi$.
For each of these we will use \cite[Equation (10) on page 240]{MacLane}
to write the relevant Kan extension as a colimit, and we will then use the 
following observation, whose proof is left to the reader.

\begin{lemma}
\label{colim}
Let $\A$ and $\B$ be categories and let $G:\A\rightarrow \mbox{Top}$
and $H:\B\rightarrow \mbox{Top}$ be functors.  Each pair consisting of a
functor 
$K:\A\rightarrow
\B$ and a natural transformation $\nu:G\rightarrow H\circ K$ induces a map
\[
\mbox{\rm colim}_\A \,\,G \rightarrow \mbox{\rm colim}_\B \,\,H
\]
\end{lemma}

We begin by constructing the transformation $\sigma_*$.  Let 
$X_1^\b,\ldots,X_k^\b$ be augmented cosimplicial spaces and
let $S$ be a totally ordered finite set;  we want to construct 
\[
\sigma_*:\Xi_k(X_1^\b,\ldots,X_k^\b)^S
\rightarrow
\Xi_k(X_{\sigma(1)}^\b,\ldots,X_{\sigma(k)}^\b)^S
\]

Let $\A_1$ be the category whose objects
are the diagrams
\[
\xymatrix@1{
\bar{k}
&
T
\ar[l]_f
\ar[r]^h
&
S
}
\]
where $T$ is a totally ordered finite set, $f$ is a map of sets and $h$ is an
ordered map; we denote such a diagram by $(f,h)$.  A morphism from $(f,h)$ to
$(f',h')$ is a commutative diagram
\[
\xymatrix@1{
\bar{k}
\ar[d]_{=}
&
T
\ar[l]_f
\ar[r]^h
\ar[d]_{g}
&
S
\ar[d]_{=}
\\
\bar{k}
&
T'
\ar[l]_{f'}
\ar[r]^{h'}
&
S
}
\]
where $g$ is an ordered map.
Let 
\[
G_1:\A_1\to \mbox{Top}
\]
be the functor which
takes $(f,h)$ to $\prod X_i^{f^{-1}(i)}$. 

By \cite[Equation (10) on page 240]{MacLane} we have
\begin{equation}
\label{F6}
\Xi_k(X_1^\b,\ldots,X_k^\b)^S= 
\mbox{colim}_{\A_1} \,\,G_1
\end{equation}
Next let $\sigma\in\Sigma_k$ and let
$H_1$ be the functor which 
takes
$(f,h)$ to $\prod X_{\sigma(i)}^{f^{-1}(i)}$;
we have
\[
\Xi_k(X_{\sigma(1)}^\b,\ldots,X_{\sigma(k)}^\b)^S=
\mbox{colim}_{\A_1} \,\,H_1
\]

\begin{definition}
\label{m1}
{\rm
The map
\[
\sigma_*:\Xi_k(X_1^\b,\ldots,X_k^\b)^S
\rightarrow
\Xi_k(X_{\sigma(1)}^\b,\ldots,X_{\sigma(k)}^\b)^S
\]
is induced by the functor
$K_1:\A_1\to\A_1$ which takes $(f,h)$ to $(\sigma^{-1}\circ f,h)$ and the 
natural transformation 
\[
\nu_1:G_1\rightarrow H_1\circ K_1
\]
which takes
$(x_1,\ldots,x_k)$ to $(x_{\sigma(1)},\ldots,x_{\sigma(k)})$.
}
\end{definition}

Next let $j_1,\ldots,j_k\geq 0$, 
let $X^\b_1,\ldots,X^\b_{j_1+\cdots+j_k}$
be augmented cosimplicial spaces and let $S$ be a finite totally ordered set.
We want to construct
\[
\Gamma_{j_1,\ldots,j_k}:
\Xi_k(\Xi_{j_1}(X_1^\b,\ldots),\ldots,\Xi_{j_k}(\ldots,X^\b_{j_1+\cdots+j_k}))^S
\rightarrow
\Xi_{j_1+\cdots+j_k}(X_1^\b,\ldots,X_{j_1+\cdots+j_k}^\b)^S
\]
First observe that 
\[
\Xi_k(\Xi_{j_1}(X_1^\b,\ldots),\ldots)^S=
\mbox{colim}_{\A_1} \,\,H,
\]
where $\A_1$ is the category defined above and
$H$ is the functor which takes $(f,h)$ to 
\[
\prod_{i=1}^k\,
\Xi_{j_i}(X_{j_1+\cdots+j_{i-1}+1}^\b,\ldots,X^\b_{j_1+\cdots+j_i})^{f^{-1}(i)}
\]
Thus a point in $\Xi_k(\Xi_{j_1}(X_1^\b,\ldots),\ldots)^S$ is
an equivalence class represented by a diagram 
\[
\xymatrix@1{
\bar{k}
&
T
\ar[l]_f
\ar[r]^h
&
S
}
\]
together with points 
$x_i\in \Xi_{j_i}(X_{j_1+\cdots+j_{i-1}+1}^\b,\ldots)^{f^{-1}(i)}$ 
for $1\leq i\leq k$. Similarly, 
each $x_i$ is represented by a diagram
\[
\xymatrix@1{
\bar{j_i}
&
T_i
\ar[l]_{f_i}
\ar[r]^-{h_i}
&
f^{-1}(i)
}
\]
together with a point 
\[
x_i'\in\prod_{p=j_1+\cdots+j_{i-1}+1}^{j_1+\cdots+j_i} X_p^{f_i^{-1}(p-j_1-\cdots-j_{i-1})}
\]
We can assemble this information into a diagram
\begin{equation}
\label{F5}
\xymatrix{
\overline{\jmath_1+\cdots+\jmath_k}
&
\bar{\jmath}_1\coprod \ldots\coprod \bar{\jmath}_k
\ar[l]_-{\chi}
\ar[d]_{\psi}
& U
\ar[l]_-e
\ar[d]_g
&
\\
&
\bar{k}
&
T
\ar[l]_f
\ar[r]^h
&
S
}
\end{equation}
and a point $x\in \prod X_p^{(\chi\circ e)^{-1}(p)}$.  Here $\chi$ is the 
unique ordered bijection, $\psi$ takes $\bar{\jmath}_i$ to $i$, $U$ is 
$\coprod T_i$, $e$ is $\coprod f_i$, the restriction of $g$ to $T_i$ is 
$h_i$, and we give $U$ the unique total order for which $g$ and the 
inclusions of the $T_i$ are ordered maps.  

Let $\A_2$ be the category of 
diagrams of the form \eqref{F5} for which $g$ and $h$ are ordered; an object 
of $\A_2$ 
will be denoted $(e,f,g,h)$.  What we have shown so far is that
$\Xi_k(\Xi_{j_1}(X_1^\b,\ldots),\ldots)^S$ is the colimit over $\A_2$
of the functor $G_2$ which takes
$(e,f,g,h)$ to $\prod X_i^{(\chi\circ e)^{-1}(i)}$.
Next we note that $\Xi_{j_1+\cdots+j_k}(X_1^\b,\ldots,X_{j_1+\cdots+j_k}^\b)^S$ is the
colimit over the category $\B$ of diagrams
\[
\xymatrix@1{
\overline{\jmath_1+\cdots+\jmath_k}
&
T
\ar[l]_-f
\ar[r]^h
&
S
}
\]
(denoted $(f,h)$) of the functor $H_2$ which takes $(f,h)$ to $\prod
X_i^{f^{-1}(i)}$.  Now let 
\[
K_2:\A_2\rightarrow\B
\]
take $(e,f,g,h)$ to 
$(\chi\circ e,h\circ g)$;  note that $H_2\circ K_2=G_2$.

\begin{definition}
\label{m2}
{\rm
The map
\[
\Gamma_{j_1,\ldots,j_k}:
\Xi_k(\Xi_{j_1}(X_1^\b,\ldots),\ldots)^S
\rightarrow
\Xi_{j_1+\cdots+j_k}(X_1^\b,\ldots,X_{j_1+\cdots+j_k}^\b)^S
\]
is induced by the functor $K_2:\A_2\to \B$ and the identity natural
transformation from $G_2$ to $H_2\circ K_2$.
}
\end{definition}

Finally,
let $K_3:\B\rightarrow \A_2$ be the functor which takes $(f,h)$ to
$(f,\psi\circ f,\mbox{id},h)$.  
Then $G_2\circ K_3=H_2$ and we can let $\nu_3:H_2\rightarrow G_2\circ K_3$ be 
the identity natural transformation; the result is a natural transformation
\[
\Lambda:
\Xi_{j_1+\cdots+j_k}(X_1^\b,\ldots,X_{j_1+\cdots+j_k}^\b)^S
\rightarrow
\Xi_k(\Xi_{j_1}(X_1^\b,\ldots),\ldots)^S
\]

With these definitions it is easy to check that $\sigma_*$, $\Gamma$ and
$\Lambda$ are
natural in $S$, that conditions (a)--(e) of definition \ref{funop} are 
satisfied, and that $\Lambda$ is inverse to $\Gamma$.  We have now shown

\begin{theorem}
\label{fo}
The collection $\Xi_k$, $k\geq 0$, with the structure maps $\sigma_*$ and
$\Gamma$ defined above, is a strict functor-operad in the category
of augmented cosimplicial spaces.  In particular $\Xi_2$ is a symmetric 
monoidal product with unit $\Xi_0$.
\qed
\end{theorem}

\begin{remark}
\label{revision1}
{\rm
If $\A$ is any symmetric monoidal category with the property that the symmetric
monoidal product preserves colimits (which is automatic when $\A$ is a closed 
symmetric monoidal category) then Theorem \ref{fo} and its proof are valid 
for the category of augmented cosimplicial objects in $\A$, with Cartesian 
products in Top replaced by the symmetric monoidal product in $\A$.
}
\end{remark}

\begin{remark}
\label{aug}
{\rm
All of the constructions in this section can be imitated for the category of 
ordinary (nonaugmented) cosimplicial spaces, provided that the maps $f$ in 
the definition of $\Q_k$ are required to be surjective.  
In this setting the product $\Xi_2$ is still associative and commutative, but
it is not unital:
$\Xi_0$ is empty, and neither $\Xi_0$ nor any other
cosimplicial space is an identity object for $\Xi_2$.
}
\end{remark}

%\newpage

\section{A sufficient condition for $\Tot(X^\b)$ to be an $E_\infty$ space.}

\label{sec4a}

If $X^\b$ is an augmented cosimplicial space, we define $\Tot(X^\b)$ to be 
the usual Tot of the cosimplicial space obtained by restricting $X^\b$ to 
$\Delta$. Equivalently, $\Tot(X^\b)$ is 
\[
\Hom(\Delta^\b,X^\b)
\]
where Hom denotes maps of augmented cosimplicial spaces and
we extend $\Delta^\b$ to an augmented cosimplicial space
by setting $\Delta^\emptyset=\emptyset$.

Now apply Proposition \ref{3.2}, letting 
$\F$ be the functor operad $\Xi$ constructed in 
Section \ref{sec4} and
$A$ the augmented cosimplicial space 
$\Delta^\b$.  This gives an operad 
$\D$ with $k$-th space $\D(k)=\Tot(\Xi_k(\Delta^\b,\ldots,\Delta^\b))$.

Recall that we have defined $\boxtimes$ to be $\Xi_2$; this is a symmetric
monoidal product on the category of augmented cosimplicial spaces.

\begin{theorem}
\label{4.5}
{\rm (a)} 
$\D$ is an $E_\infty$ operad.

{\rm (b)} 
If $X^\b$ is a commutative monoid with respect to $\boxtimes$
(equivalently, if $X^\b$ is an algebra over $\Xi$) then $\D$ acts on 
Tot$(X^\b)$.  
\end{theorem}

\Proof Part (b) is immediate from Propositions \ref{3.4} and \ref{3.6}.  Part
(a) will be proved in Section \ref{secA}. \QED

\begin{remark}
{\rm
The analog of Theorem \ref{4.5}(b) is valid, with the same proof, when
$X^\b$ is an augmented cosimplicial spectrum. 
}
\end{remark}

In the remainder of this section we use the adjointness property of
$\lan_\Phi$
\cite[beginning of Section X.3]{MacLane} to give an explicit characterization
of commutative $\boxtimes$-monoids.  

\begin{definition}
{\rm
Let $X^\b$ be an augmented cosimplicial space.  A 
$\langle \  \rangle$-structure 
on $X^\b$ 
consists of a map
\[
\langle f \rangle:X^{f^{-1}(1)}\times X^{f^{-1}(2)}\to X^T
\]
for each totally ordered set $T$ and each $f:T\to \{1,2\}$.
}
\end{definition}

\begin{definition}
{\rm
A 
$\langle \  \rangle$-structure
on $X^\b$ is {\it consistent} if
for every commutative diagram
\[
\xymatrix{
T
\ar[rr]^h
\ar[dr]_-f
&&
T'
\ar[dl]^-g
\\
&
\{1,2\}
&
}
\]
the diagram 
\[
\xymatrix{
X^{f^{-1}(1)}\times X^{f^{-1}(2)}
\ar[d]_{(h_1)_*\times (h_2)_*}
\ar[r]^-{\langle f\rangle}
&
X^T 
\ar[d]^{h_*}
\\
X^{g^{-1}(1)}\times X^{g^{-1}(2)}
\ar[r]^-{\langle g\rangle}
&
X^{T'}
}
\]
commutes, where $h_i$ is the restriction of $h$ to $f^{-1}(i)$.
}
\end{definition}

Recall the notation of Definition \ref{Kan2}.
Since $\lan_\Phi$ is left adjoint to $\Phi^*$, a map $X^\b\boxtimes X^\b\to
X^\b$ is the same thing as a natural transformation
\[
(X^\b\bar\times X^\b)\circ \Psi \to X^\b \circ \Phi
\]
and it's easy to check that this is the same thing
as a consistent $\langle\ \rangle$-structure on $X^\b$.  It remains to 
translate the commutativity, associativity and unitality conditions satisfied
by a commutative $\boxtimes$-monoid into this language.

\begin{definition}
{\rm
A 
$\langle \  \rangle$-structure
on $X^\b$ is {\it commutative} if
the diagram 
\[
\xymatrix{
X^{f^{-1}(1)}\times X^{f^{-1}(2)}
\ar[d]_{\tau}
\ar[r]^-{\langle f\rangle}
&
X^T 
\ar[d]^{=}
\\
X^{f^{-1}(2)}\times X^{f^{-1}(1)}
\ar[r]^-{\langle t\circ f\rangle}
&
X^T
}
\]
commutes, where $\tau$ is the switch map and $t$ is the transposition of
$\{1,2\}$.
}
\end{definition}

For the associativity condition we need some notation.  Let $T$ be a totally 
ordered set and let $g:T\to \{1,2,3\}$ be a function.  Define
\[
\alpha:\{1,2,3\}\to\{1,2\} 
\]
by $\alpha(1)=1,\alpha(2)=1,\alpha(3)=2$ and define
\[
\beta:\{1,2,3\}\to\{1,2\}
\]
by $\beta(1)=1,\beta(2)=2,\beta(3)=2$.
Let $g_1$ be the restriction of $g$ to $g^{-1}\{1,2\}$ and let $g_2$ be the
restriction of $g$ to $g^{-1}\{2,3\}$.

\begin{definition}
{\rm
A 
$\langle \  \rangle$-structure
on $X^\b$ is {\it associative} if, with the
notation above, the diagram 
\[
\xymatrix{
X^{g^{-1}(1)} \times X^{g^{-1}(2)} \times X^{g^{-1}(3)} 
\ar[d]_{1\times \langle g_2\rangle}
\ar[r]^-{\langle g_1\rangle\times 1}
&
X^{g^{-1}\{1,2\}} \times X^{g^{-1}(3)}
\ar[d]^{\langle \alpha \circ g\rangle}
\\
X^{g^{-1}(1)} \times X^{g^{-1}\{2,3\}}
\ar[r]^-{\langle \beta\circ g\rangle}
&
X^T
}
\]
commutes for every choice of $T$ and of $g:T\to \{1,2,3\}$.
}
\end{definition}

\begin{definition}
{\rm
A 
$\langle \  \rangle$-structure
on $X^\b$ is {\it unital} if
there is an element $\varepsilon\in X^\emptyset$ with the property that if 
$f:T\to\{1,2\}$ takes all of $T$ to 1 then $\langle f\rangle (x,\varepsilon)=x$
for all $x$ and if $f$ takes all of $T$ to 2 then $\langle f\rangle
(\varepsilon,x)=x$ for all $x$.
}
\end{definition}

\begin{proposition}
\label{MM1}
A commutative $\boxtimes$-monoid structure on $X^\b$ determines, and is
determined by, a $\langle\ \rangle$-structure on $X^\b$ which is
consistent, commutative, associative and unital.
\end{proposition}

The proof is a routine verification using the definitions in Section
\ref{sec4}.

\begin{remark}
{\rm
The $\langle f\rangle$ operations on $S^\b W$ defined in Section \ref{sec3a}
give a consistent, commutative, associative and unital $\langle\
\rangle$-structure on $S^\b W$.
}
\end{remark}

%\newpage

\section{A filtration of $\Xi$ by functor-operads.}

\label{sec5}

In this section we describe a filtration of $\Xi$ by functor-operads
$\Xi^n$; the operad associated to $\Xi^n$ will turn out to be
equivalent to the little $n$-cubes operad $\C_n$.  

We begin with some motivation.  If $T$ is a totally ordered set and $f:T\to 
\{1,2\}$ is a function, the two totally ordered sets $f^{-1}(1)$ and 
$f^{-1}(2)$ are mixed together to form $T$.  The amount of mixing can be 
measured by the number of times the value of $f$ switches from 1 to 2 or from 
2 to 1 as one moves through the set $T$.  The idea in the definition of 
$\Xi^n$ is to control the amount of mixing that is allowed.  

\begin{definition}
\label{complexity}
{\rm
Let $T$ be a finite totally ordered set, let $k\geq 2$, and let
$f:T\rightarrow \bar{k}$.
We define the {\it complexity} of $f$ as follows.
If $k$ is $0$ or $1$ the complexity is 0.
If $k=2$ let $\sim$ be
the equivalence relation on $T$ generated by
\[
a\sim b \mbox{ if $a$ is adjacent to $b$ and $f(a)=f(b)$}
\]
and define the complexity of $f$ to be the number of equivalence classes
minus 1.  If $k> 2$  define the complexity of $f$ to be the maximum of the
complexities of the restrictions $f|_{f^{-1}(A)}$ as $A$ ranges over the
two-element subsets of $\bar{k}$.
}
\end{definition}

\begin{remark}
{\rm
If $k=2$ the complexity of $f$ is exactly the amount of mixing in $f$ as
discussed above.
The definition of complexity is suggested by \cite{Smith};
the reason we use it here is that it is well-adapted to the proofs of Theorems
\ref{S2} and \ref{6.4}(a).
There may be other ways of defining complexity that would also lead to 
Theorems \ref{S2} and \ref{6.4}(a), although this seems unlikely.
}
\end{remark}

Now fix $n\geq 1$. Recall the category 
$\Q_k$ from Definition \ref{Q}.

\begin{definition}
\label{MM3}
{\rm
Let $\Q_k^n$ be the full 
subcategory of $\Q_k$ whose objects are pairs $(f,S)$ where $f$ has 
complexity $\leq n$.  Let
\[
\iota^n:\Q_k^n \to \Q_k
\]
be the inclusion.
}
\end{definition}

\begin{definition}
\label{Kan3}
{\rm
For each $n\geq 1$ and each $k\geq 0$ define a functor $\Xi_k^n$ as follows.
Given augmented cosimplicial spaces $X_1^\b,\ldots,X_k^\b$,
let $X_1^\b\bar\times\cdots\bar\times X_k^\b$ be the functor defined in
Definition \ref{Kan2} and let
$\Xi_k^n(X_1^\b,\ldots,X_k^\b)$ 
be the Kan extension
\[
\mbox{Lan}_{\Phi\circ\iota^n} ((X_1^\b\bar\times\cdots\bar\times 
X_k^\b)\circ \Psi\circ\iota^n)
\]
}
\end{definition}

\begin{theorem}
\label{S2}
For each $n\geq 1$, $\Xi^n$ is a (non-strict) functor-operad in the category 
of augmented cosimplicial spaces.
\end{theorem}

\Proof  
We define $\sigma_*$ as in the proof of Theorem \ref{fo}, using the fact 
that the complexity of $\sigma^{-1}\circ f$ is the same as that of $f$. 
We define $\Gamma_{j_1,\ldots,j_k}$ as in the proof of Theorem \ref{fo},
but we must verify that if the complexities of
$e|_{e^{-1}(\bar{\jmath}_i)}$ and $f$ in diagram \eqref{F5}
are $\leq n$ then the complexity of $\chi\circ e$ will be also be $\leq n$.
For this we need to show that for each two-element subset $A$ of
$\bar{\jmath}_1\coprod\ldots\coprod\bar{\jmath}_k$ the complexity
of $\chi\circ e|_{e^{-1}(A)}$ will be $\leq n$; but this is true when $A$ is 
contained in some $\bar{\jmath}_i$ (because the complexity of
$e|_{e^{-1}(\bar{\jmath}_i)}$ is $\leq n$) and it is also true if $A$ is not
contained in any $\bar{\jmath}_i$ (because the complexity of $f$ is
$\leq n$).
\QED

\begin{remark}
{\rm
If $\A$ is any symmetric monoidal category satisfying the hypothesis of Remark
\ref{revision1} then
Theorem \ref{S2} and its proof are valid for the category of augmented 
cosimplicial objects in $\A$, with Cartesian 
products in Top replaced by the symmetric monoidal product in $\A$.
}
\end{remark}

\begin{remark}
\label{MM2}
{\rm
In the special case $n=1$, the functor-operad $\Xi^1$ is closely related to the
monoidal product $\bx$ defined in Section \ref{sec2}.  First observe that
order-preserving maps $f:T\to\bar{k}$ have filtration 1 and that every map of
filtration 1 can be written uniquely as the composite of an order-preserving 
map and a permuation of $\bar{k}$.  If we use
order-preserving maps in Definitions \ref{MM3} and  \ref{Kan3} instead of 
maps of filtration 1, we get a nonsymmetric strict functor-operad $\Upsilon$ 
which is related to both $\bx$ and $\Xi^1$:

(a) The restriction of
$\Upsilon_k$ to the category of (unaugmented) cosimplicial spaces 
is naturally isomorphic to the iterated $\bx$-product $\bx^k$.

(b)
$\Xi^1_k$ is naturally isomorphic to
\[
\coprod_{\sigma\in \Sigma_k} \Upsilon_k \circ \sigma_\# 
\]
That is, $\Xi^1$ is obtained by extending the nonsymmetric functor-operad 
$\Upsilon$ in the obvious way to a (symmetric) functor-operad.
}
\end{remark}

%\newpage

\section{An operad which acts on Tot of a $\Xi^n$-algebra.}

\label{sec5a}

Applying Proposition \ref{3.2} with $\F=\Xi^n$ and $A=\Delta^\b$ we get an 
operad $\D_n$ with $k$-th space 
$\D_n(k)=\Tot(\Xi^n_k(\Delta^\b,\ldots,\Delta^\b))$.

\begin{theorem}
\label{6.4}
{\rm (a)} 
$\D_n$ is weakly equivalent in the category of operads to 
$\C_n$.

{\rm (b)} 
If $X^\b$ is an algebra over $\Xi^n$ then $\D_n$ acts on
Tot$(X^\b)$.
\end{theorem}

The statement of part (a) means that there is a chain of operads and weak
equivalences of operads
\[
\D_n \leftarrow \cdots \rightarrow \C_n
\]
Part (b) of the Theorem is immediate from Propositions \ref{3.4} and 
\ref{3.6}. Part (a) will be proved in Section \ref{secJ}.

\begin{remark}
{\rm
The analog of Theorem \ref{6.4}(b) is valid, with the same proof, when
$X^\b$ is an augmented cosimplicial spectrum. 
}
\end{remark}

In the remainder of this section we give an explicit characterization of 
$\Xi^n$-algebras, analogous to that given in Section \ref{sec4a} for
commutative $\boxtimes$-monoids.  

\begin{definition}
{\rm
Let $X^\b$ be an augmented cosimplicial space.  An $n$-structure
on $X^\b$
consists of a map
\[
\langle f \rangle:X^{f^{-1}(1)}\times \cdots \times
X^{f^{-1}(k)}
\to X^T
\]
for each totally ordered set $T$, each $k\geq 0$, and each $f:T\to \bar{k}$ 
with complexity
$\leq n$.
}
\end{definition}

\begin{definition}
\label{new2}
{\rm
An $n$-structure
on $X^\b$ is {\it consistent} if,
for every commutative diagram
\[
\xymatrix{
T
\ar[rr]^h
\ar[dr]_-f
&&
T'
\ar[dl]^-g
\\
&
\bar{k}
&
}
\]
in which $f$ and $g$ have complexity $\leq n$,
the diagram 
\[
\xymatrix{
\prod_{i=1}^k X^{f^{-1}(i)}
\ar[d]_{\prod (h_i)_*}
\ar[r]^-{\langle f\rangle}
&
X^T 
\ar[d]^{h_*}
\\
\prod_{i=1}^k X^{g^{-1}(i)}
\ar[r]^-{\langle g\rangle}
&
X^{T'}
}
\]
commutes, where $h_i$ is the restriction of $h$ to $f^{-1}(i)$.
}
\end{definition}

It's easy to check (using the fact that $\lan_\Phi$ is left adjoint to
$\Phi^*$) that a consistent $n$-structure on $X^\b$ is the same thing as a 
collection of maps
\[
\Xi^n_k(X^\b,\ldots,X^\b)\to X^\b,
\]  
one for each $k\geq 0$.
It remains to
translate the rest of the definition of $\Xi^n$-algebra into this language.

\begin{definition}
\label{new3}
{\rm
An $n$-structure
on $X^\b$ is {\it commutative} if, for each $f$ with complexity $\leq n$ and
each $\sigma\in \Sigma_k$,
the diagram
\[
\xymatrix{
\prod_{i=1}^k X^{f^{-1}(i)}
\ar[d]_{s}
\ar[rr]^-{\langle f\rangle}
&&
X^T
\ar[d]^{=}
\\
\prod_{i=1}^k X^{f^{-1}(\sigma(i))}
\ar[rr]^-{\langle \sigma^{-1}\circ f\rangle}
&&
X^T
}
\]
commutes  (where $s$ is the evident
permutation of the factors). 
}
\end{definition}

For the next definition we need some notation. 
Suppose we are given a partially ordered set $T$, numbers 
$k,j_1,\ldots,j_k\geq 0$, and maps
\[
f:T\to \bar{k}
\]
and 
\[
g_i:f^{-1}(i)\to \ov{\jmath_i}
\]
for $1\leq i\leq k$.  Let $j=\sum j_i$.  The maps $g_i$ determine a map
\[
g:T\to \bar{\jmath}
\]
in an evident way; the formula for $g$ is 
\[
g(a)=g_i(a) +\sum_{i'<i} j_{i'} \quad\mbox{if $a\in f^{-1}(i)$}
\]

\begin{definition}
\label{new4}
{\rm
An $n$-structure on $X^\b$ is {\it associative} if
the following diagram commutes for every choice of $f$ and $g_1,\ldots,g_k$
with complexity $\leq n$:
\[
\xymatrix{
\prod_{i=1}^k \prod_{b=1}^{j_i} X^{g_i^{-1}(b)}  \
\ar[d]_{=}
\ar[r]^-{\prod \langle g_i \rangle}
&
\ \prod_{i=1}^k X^{f^{-1}(i)} 
\ar[d]^{\langle f\rangle}
\\
\prod_{c=1}^j X^{g^{-1}(c)} 
\ar[r]^-{\langle g\rangle}
&
X^T 
}
\]
}
\end{definition}

In order to state the unitality condition we need some more notation.  If $i\in
\bar{k}$ let $\lambda_i:\overline{k-1}\to \bar{k}$ be the order-preserving monomorphism whose image does not contain $i$.

\begin{definition}
\label{new5}
{\rm
An $n$-structure
on $X^\b$ is {\it unital} if
there is an element $\varepsilon\in X^\emptyset$ with the following property: 
\[
\langle \lambda_i\circ f\rangle 
(x_1,\ldots,x_{i-1},\varepsilon,x_{i},\ldots,x_{k-1})
=\langle f\rangle (x_1,\ldots,x_{i-1},x_{i},\ldots,x_{k-1})
\]
for all $f:T\to \overline{k-1}$ with complexity $\leq n$, all $i\in \bar{k}$, 
and all choices of 
$x_1,\ldots,x_{k-1}$.
}
\end{definition}

\begin{proposition}
\label{new1}
A $\Xi^n$-algebra structure on $X^\b$ determines, and is
determined by, an $n$-structure on $X^\b$ which is
consistent, commutative, associative and unital.
\end{proposition}

The proof is a routine verification using the definitions in Section
\ref{sec4}.

%\newpage

\section{Example: the cosimplicial space associated to a nonsymmetric
operad with multiplication.}

\label{sec6}

Let us say that an augmented
cosimplicial space $X^\b$ is {\it reduced} if $X^\emptyset$ is a point.

In this section we specialize to the case $n=2$. 
We show that a $\Xi^2$ structure on a reduced augmented cosimplicial 
space is the same thing as
a ``nonsymmetric operad with multiplication.'' 

First recall \cite[Definition II.1.14]{MSS} that the definition of {\it
nonsymmetric operad} is obtained from the usual definition of operad by
deleting all references to symmetric groups.  

Let {\it Ass} be the nonsymmetric operad whose $k$-th space is a point for
all $k\geq 0$.  

The next definition is due to Gerstenhaber and Voronov \cite{GV}.

\begin{definition}
\label{om}
{\rm
A {\it nonsymmetric operad with multiplication} is a nonsymmetric operad 
$\O$ together with a morphism $\text{\it Ass}\to \O$. 
}
\end{definition}

\begin{remark}
{\rm
i) It is easy to check that a morphism $\text{\it Ass}\to \O$ is the same
thing as a pair of elements $\mu\in\O(2)$, $e\in O(0)$ satisfying 
\[
\gamma(\mu,\mu,\id)=\gamma(\mu,\id,\mu)
\]
and
\[
\gamma(\mu,e,\id)=\gamma(\mu,\id,e)=\id
\]
where $\gamma$ is the composition operation and $\id$ is the identity element
of the nonsymmetric operad $\O$.

ii) An important example of a nonsymmetric operad with multiplication (in the 
category of abelian groups) is
\[
\O(k)=\Hom_{\mathbb Z}(A^{\otimes k}, A)
\]
where $A$ is an associative ring.  Here the composition operation is the 
obvious one, $\id\in\O(1)$ is the identity map of the ring $A$, $\mu\in 
\O(2)$ is the multiplication of the ring $A$, and $e\in \O(0)$ is the 
identity element of the ring $A$.  This example is related to the Hochschild 
cochain complex. 
}
\end{remark}

The rest of this section is devoted to the proof of:

\begin{proposition}
A $\Xi^2$ structure on a reduced augmented 
cosimplicial space $X^\b$ determines, and is determined by, a structure of
nonsymmetric operad with multiplication on the sequence of spaces $X^k$, 
$k\geq 0$.
\end{proposition}

First suppose that $X^\b$ is a reduced augmented $\Xi^2$-algebra; we need to
define the composition operation $\gamma$, the identity element $\id\in X^1$, 
and elements $\mu\in X^2$ and $e\in X^0$.

For the composition operation, first note that
if $U_1,U_2,\ldots,U_m$ are finite totally ordered sets there is a unique
total order on $U_1\,\coprod \,\cdots \,\coprod 
\,U_m$ 
for which the inclusion maps into the coproduct are order-preserving and 
every element of $U_i$ is less 
than every element of $U_j$ for $i<j$.  Now let $k,j_1,\ldots,j_k\geq 0$ be 
given, and let $T$ be 
the totally ordered set
\[
\{0\} \,\coprod \,[j_1] \,\coprod
\,\{1\} \,\coprod \,[j_2] \,\coprod
\,\ldots\, \coprod
\,[j_k]\,\coprod \,\{k\}
\]
Let $f:T\to \overline{k+1}$ be the map that takes each set $\{i\}$ to 1 and 
each $[j_i]$ to $i+1$; thus $f^{-1}(1)=[k]$ and $f^{-1}(i)=[j_{i-1}]$ for
$i>1$.
Let $\sim$ be the equivalence relation on $T$ 
generated 
by: $x\sim y$ if $x$ is adjacent to $y$ in the total order and
$f(x)\neq f(y)$.  Let $S$ be the quotient $T/\sim$.
$S$ inherits a total order from $T$ and has $j_1+\cdots+j_k+1$ elements.
Let $h$ be the composite
\[
T\to S\to [j_1+\cdots+j_k],
\]
where the first map is the quotient map and the second is the unique
order-preserving bijection.
Finally, let
\[
\gamma:X^k\times X^{j_1} \times \cdots X^{j_k}\to X^{j_1+\cdots+j_k}
\]
be the composite $h_*\circ\langle f\rangle$.  

Next let $\varepsilon$ be the unique element of $X^\emptyset$ and define the 
elements $e$, $\id$ and $\mu$ by $e=d^0 \varepsilon$, $\id=d^0 e$, and
$\mu=d^0\id$.

It is easy to check that these definitions give the sequence of spaces $X^k,
k\geq 0$ a structure of nonsymmetric operad with multiplication.

In the other direction, suppose that $\O$ is a nonsymmetric operad with 
multiplication.  We begin by recalling a standard piece of notation:
if $x\in\O(k), y\in\O(j)$ and $1\leq i\leq k$, define
\[
x\circ_i y=\gamma(x,\id,\ldots,y,\ldots\id)
\]
where the $y$ is preceded by $i-1$ copies of id.  This gives a map
\[
\circ_i:\O(k)\times \O(j)\to \O(k+j-1)
\]
which inherits associativity and unital properties from $\gamma$ (see 
\cite[page 46]{MSS} for details). 

We can now define the augmented cosimplicial space $\O^\b$ associated to $\O$.
Let $\O^\emptyset$ be a single point $\varepsilon$. 
Define the cosimplicial structure maps by letting
$d^0\varepsilon = e$ and, if $x\in \O^p$ with $p\geq 0$,
\begin{align*}
d^i x &=
\left\{
\begin{array}{ll}
\mu\circ_2 x & \mbox{if $i=0$} \\
x\circ_i \mu & \mbox{if $0<i<p+1$} \\
\mu\circ_1 x & \mbox{if $i=p+1$}
\end{array}
\right. \\
s^i x  &= x\circ_{i+1} e. 
\end{align*}
(This definition is motivated by the definition of the structure maps in the
Hochschild cochain complex.)

It remains to define $\langle f \rangle$ operations on $\O^\b$.  

For this we need some preliminary definitions.  Define
\[
\sqcup:\O^p\times\O^q\to\O^{p+q+1}
\]
by
\[
x\sqcup y=\gamma(\mu,x,d^0y)
\]

\begin{remark}
{\rm
A function $f:T\to \bar{k}$ is the same thing as a finite sequence with values
in $\bar{k}$.  If
$f:[p+q+1]\to\bar{2}$ is the
sequence $1\ldots 12\ldots 2$ with 1 repeated $p+1$ times
and 2 repeated $q+1$ times
then $\langle f\rangle$ will be defined to be $\sqcup$.
}
\end{remark}

Next define
\[
\tilde{\gamma}: \O(k)\times\O(j_1)\times\cdots\times\O(j_k)
\to
\O(2k+j_1+\cdots+j_k)
\]
by
\[
\tilde{\gamma}(x,y_1,\ldots,y_k)
= \gamma(x,d^0d^{j_1+1}y_1,\ldots,d^0d^{j_k+1}y_k)
\]

\begin{remark}
{\rm
If $f:[2k+j_1+\cdots+j_k]\to \overline{k+1}$ is the sequence 
\[
12\cdots 213\cdots 31\cdots 1\ k+1\cdots \ k+1\ 1,
\]
with each $i>1$ repeated $j_{i-1}$ times,
then $\langle f\rangle$ will be defined to be $\tilde{\gamma}$.
}
\end{remark}

Next note that if $S$ is a finite totally ordered set there is a 
canonical homeomorphism $\O^S\cong\O^{||S||}$, where $||S||$ is the number of 
elements in $S$ minus 1.

Now let 
\[
f:T\to \bar{k}
\]
be a map of complexity $\leq 2$.  We will define $\langle f \rangle$ by
induction on $||T||$.

First of all, if $T$ is empty, then $\langle f\rangle$ is defined to be the
unique map
\[
\O^\emptyset\times\cdots\times\O^\emptyset\to \O^\emptyset
\]
If $T$ has a single element then $\langle f \rangle$ is the canonical
homeomorphism
\[
\O^\emptyset\times\cdots\times O^T\times\cdots\cong O^T
\]

Next define a {\it segment} of $f$ to be a subset $S$ of $T$ such that $f$ 
has the same value on the minimal and maximal elements of $S$, and define a 
{\it maximal} segment to be a segment which is not properly contained in any 
other segment.  Let $S_1,\ldots,S_r$ denote the maximal segments of $f$; then
$T$ is the union of the $S_j$, and the fact that $f$ has complexity $\leq 2$ 
implies that the $S_j$ are disjoint.  Also note that each $f^{-1}(i)$ is
contained in some $S_j$.

If $r>1$ (that is, if $f$ has more than one maximal segment), let 
$g_1,\ldots,g_r$ be the restrictions of $f$ to $S_1,\ldots,S_r$ and define 
$\langle f \rangle$ to be the composite
\[
\prod_{i=1}^k \O^{f^{-1}(i)}
\cong
\prod_{j=1}^r\prod_{i\in g_j(S_j)} \O^{g_j^{-1}(i)}
\labarrow{\prod_{j=1}^r \langle g_j\rangle}
\O^{S_1}\times\cdots\times\O^{S_r}
\labarrow{\sqcup}
\O^T
\]

If $r=1$ let $j$ be the value of $f$ at the minimum and maximum elements of $T$
and let $t_0,\ldots,t_{||f^{-1}(j)||}$ be the elements of $f^{-1}(j)$ in 
increasing order.  For each $l$ from 1 to $||f^{-1}(j)||$ let $U_l$ be the set
$\{ t_{l-1}<t<t_l \}$; then $T$ is the disjoint union of $f^{-1}(j)$ and the
sets $U_l$, and each $f^{-1}(i)$ is contained in one of the pieces of this
disjoint union.  Let $g_l$ be the restriction of $f$ to $U_l$.
We define $\langle f \rangle$ to be the composite
\[
\prod_{i=1}^k \O^{f^{-1}(i)}
\cong
\O^{f^{-1}(j)}\times \prod_l \prod_{i\in g_l(U_l)} \O^{g_l^{-1}(i)}
\labarrow{1\times \prod_l \langle g_l \rangle}
\O^{f^{-1}(j)}\times \prod_l \O^{U_l}
\labarrow{\tilde{\gamma}}
\O^T 
\]

Now it is straightforward, although tedious, to check that the 
operations $\langle f \rangle$ that we have defined satisfy the hypothesis of 
Proposition \ref{new1} with $n=2$.

\section{Example: $\Omega^n Y$}

\label{sec7}

One of the basic properties of the little $n$-cubes operad $\C_n$ is that it
acts naturally on $n$-fold loop spaces.  In this section we show that the
operad $\D_n$ has the same property:

\begin{proposition}
$\D_n$ acts naturally on $n$-fold loop spaces.
\end{proposition}

First we need to give a cosimplicial model for the $n$-th loop space of a
pointed space $Y$.  Let $S^n_\b$ be
the quotient of $\Delta^n_\b$ by its $(n-1)$-skeleton; then $S^n_\b$ is a
pointed simplicial set whose realization is the
$n$-sphere.  Let $L_Y^\b$ be the cosimplicial space 
$\Map_*(S^n_\b,Y)$ (when $n=1$, this is just the usual geometric cobar
construction on $Y$).

\begin{lemma}
$\Tot(L_Y^\b)$ is homeomorphic to $\Omega^n Y$.
\end{lemma}

\Proof
This is easy from the definition of $\Tot$, using the fact that
$S^n_\b$ has only one non-degenerate simplex other than the basepoint.
\QED

It only remains to show that $L_Y^\b$ is a $\Xi^n$-algebra.  First we observe
that any simplicial set $Z_\b$ has co-operations
\[
\lfloor f\rfloor:Z_T\to Z_{f^{-1}(1)}\times\cdots\times Z_{f^{-1}(k)}
\]
for all $f: T\to \bar{k}$, given by 
\[
\lfloor f\rfloor(x)=(h_1^*(x),\ldots,h_k^*(x))
\]
where $h_i:f^{-1}(i)\to T$ is the inclusion.
It is easy to check that these co-operations have properties dual to
Definitions \ref{new2}, \ref{new3}, \ref{new4} and \ref{new5}.

\begin{lemma}
\label{new6}
If $f$ has complexity $\leq n$ then the map
\[
\lfloor f\rfloor:S^n_T\to S^n_{f^{-1}(1)}\times\cdots\times S^n_{f^{-1}(k)}
\]
factors through the wedge
\[
S^n_{f^{-1}(1)}\vee\cdots\vee S^n_{f^{-1}(k)}
\]
\end{lemma}

\Proof
First observe that $f$ can be thought of as a sequence of 
length $|T|$ with values in $\bar{k}$, and that (because $f$ has complexity 
$\leq n$) this sequence has no subsequence of length $n+2$ of the form 
$ijij\cdots$. 

Next recall that an element of $\Delta^n_T$ is a function $\phi:T\to
[n]$, and that this element is in the $(n-1)$-skeleton if and only if $\phi$
is not onto.   Let $\phi:T\to [n]$ be onto.  The entries of 
\[
\lfloor f\rfloor(\phi)\in \Delta^n_{f^{-1}(1)}\times\cdots\times 
\Delta^n_{f^{-1}(k)}
\]
are the restrictions $\phi|_{f^{-1}(i)}$.  If two of these restrictions (say
when $i=a$ and $i=b$) were onto, then the sequence corresponding to $f$ would
have a subsequence $aba\cdots$ or $bab\cdots$ of length $n+2$, which is
impossible.  Thus all but one of the entries of $\lfloor f\rfloor(\phi)$ must
be in the $(n-1)$-skeleton of $\Delta^n_\b$, which proves the lemma.
\QED

Using Lemma \ref{new6}, we see that any $f$ of complexity $\leq n$ induces a map
\[
\langle f\rangle:
\prod L_Y^{f^{-1}(i)}
=
\Map_*(\bigvee S^n_{f^{-1}(i)},Y)
\to
\Map_*(S^n_T,Y) 
=
L_Y^T 
\]
These maps satisfy Definitions \ref{new2}, \ref{new3}, \ref{new4} and
\ref{new5} because the $\lfloor f\rfloor$'s satisfy the duals of those
definitions.  This completes the proof of the Proposition.

\section{The structure of $\Xi_k(\Delta^\b,\ldots,\Delta^\b)$.}
\label{secA}

Throughout this section we write $Y_k^\b$ for the augmented cosimplicial 
space 
$\Xi_k(\Delta^\b,\ldots,\Delta^\b)$ and $Y_k^S$ for the value of $Y_k^\b$ at the
finite totally ordered set $S$.  We want to investigate the structure of
$Y_k^S$ and $Y_k^\b$.

Recall (equation \eqref{F6}) that
$Y_k^S$ is colim$_{\A_1}\,\, G_1$, where $\A_1$ is the category defined just
before equation \eqref{F6} and
$G_1$ is the functor which takes the diagram
\[
\xymatrix{
\bar{k}
&
T
\ar[l]_{f}
\ar[r]^h
&
S
}
\]
to $\prod_{1\leq i\leq k}\, \Delta^{f^{-1}(i)}$.  

\begin{notation}
\label{nota}
{\rm
A diagram of the form 
\[
\xymatrix@1{
\bar{k}
&
T
\ar[l]_f
\ar[r]^h
&
S,
}
\]
where $h$ is ordered, will be denoted from now on by $(f,T,h)$.
}
\end{notation}

The elements of $Y_k^S$ are equivalence classes of pairs
$((f,T,h),u)$ 
with $u\in\prod_{i}\,\,\Delta^{f^{-1}(i)}$;  
we think of $u$ as a tuple indexed by $T$, subject to the condition that
$\sum_{a\in{f^{-1}}(i)} u_a =1$ for each $i\in\bar{k}$.
Note that $f$ must be surjective because $\Delta^\emptyset=\emptyset$.

\begin{definition}
\label{nondeg}
{\rm
(a)
A diagram
$(f,T,h)$ is
{\it nondegenerate} if $T=[m]$ for some $m$, $f$ is surjective, and for each 
$j<m$
either $f(j)\neq f(j+1)$ or $h(j)\neq h(j+1)$.

(b)
A pair $((f,T,h),u)$ with $u\in\prod_{i}\,\,\Delta^{f^{-1}(i)}$ is 
nondegenerate if $(f,T,h)$ is
nondegenerate and $u_a\neq 0$ for all $a\in T$. 
}
\end{definition}

\begin{proposition}
\label{5.1}
Each point in $Y_k^S$ is represented by a unique nondegenerate pair.
\end{proposition}

\begin{corollary}
\label{5.2}
$Y_k^S$ is a CW complex with one cell of dimension $m+1-k$ for each 
nondegenerate $(f,[m],h)$; 
the characteristic map of the cell corresponding to $(f,[m],h)$ is a
homeomorphism from $\prod\Delta^{f^{-1}(i)}$ to the closure of the cell (and
thus $Y_k^S$ is a regular CW complex). \qed
\end{corollary}

\noindent
{\bf Proof of \ref{5.1}.}
We define a function $\Upsilon$ from pairs to pairs as follows.
Given a pair $((f,T,h),u)$,
let $\sim$ be the equivalence relation on $T$ generated by 
\[
a\sim b \quad \mbox{if $a$ is adjacent to $b$, } f(a)=f(b), \mbox{ and }
h(a)=h(b)
\]
and let $T_1$ be the subset 
\[
\{ a\in T \,|\, u_a\neq 0 \}.
\]
Let $\hat{T}$ be $T_1/\sim$.  Then $f$ and $h$ induce maps
$\hat{f}:\hat{T}\rightarrow \bar{k}$ and $\hat{h}:\hat{T}\rightarrow S$.
Also, let $\pi:T_1\rightarrow \hat{T}$ be the projection and define 
$\hat{u}\in\prod_i\,\, \Delta^{\hat{f}^{-1}(i)}$ by 
$\hat{u}_c=\sum_{a\in\pi^{-1}(c)} u_a$ for $c\in \hat{T}$.
Finally, let $m=|\hat{T}|-1$ and let $g:[m]\rightarrow \hat{T}$ be the unique 
ordered
bijection.
We define $\Upsilon((f,T,h),u)=((\hat{f}\circ g,[m],\hat{h}\circ 
g),\hat{u}\circ g)$. The 
proposition is immediate from the following properties of $\Upsilon$:

(i) $\Upsilon((f,T,h),u)$ is nondegenerate.

(ii) $((f,T,h),u)$ and $\Upsilon((f,T,h),u)$ represent the same point in
$Y_k^S$.

(iii) If $((f,T,h),u)$ and $((f',T',h'),u')$ represent the same point in 
$Y_k^S$ then $\Upsilon((f,T,h),u)=\Upsilon((f',T',h'),u')$

\QED

Our next goal is to show for each $k$ that $Y_k^\b$ is isomorphic as an 
augmented cosimplicial space to $\Delta^\b\times Y_k^0$ (this is the analog for
$\boxtimes$ of Lemma \ref{iso}).

\begin{definition}
\label{eta}
{\rm
For each $S$, $\eta_S$ is the unique map 
$S\rightarrow [0]$.  This will also be denoted by $\eta$ when $S$ is clear from
the context.
}
\end{definition}

\begin{definition}
\label{omega}
{\rm
Define
\[
\omega^S:Y_k^S\rightarrow \Delta^S\times Y_k^0
\]
by letting the projection on the second factor be the map $\eta_*$ induced 
by $\eta$
and letting the projection on the first factor take
the equivalence class of $((f,T,h),u)$ to 
the element $v\in\Delta^S$
with $v_a=\frac{1}{k}\sum_{b\in h^{-1}(a)} u_b$.
The $\omega^S$ fit together to give a cosimplicial map
\[
\omega:Y_k^\b\rightarrow \Delta^\b\times Y_k^0
\]
}
\end{definition}

\begin{proposition}
\label{5.3}
$\omega$ is an isomorphism of augmented cosimplicial spaces.
\end{proposition}

\Proof
The diagram
\[
\xymatrix{
Y_k^S
\ar[dr]_{\eta_*}
\ar[rr]^{\omega^S}
&&
\Delta^S\times Y_k^0
\ar[dl]^{\pi_2}
\\
&Y_k^0&
}
\]
commutes, where $\pi_2$ is the projection.  We begin by showing 
\begin{center}
(1)\hfill for each point $y\in Y_k^0$ 
the map
$
\eta_*^{-1}(y)\rightarrow \pi_2^{-1}(y)
$
induced by $\omega^S$ is a bijection. \hfill \mbox{}
\end{center}
For this, it suffices to show 
\begin{center}
(2)\hfill the composite
$
\eta_*^{-1}(y)\labarrow{\omega^S} \Delta^S \times Y_k^0\labarrow{\pi_1} 
\Delta^S
$
is a bijection.\hfill\mbox{}
\end{center}

So let $y$ be a point of $Y_k^0$ and let
\[
((f,[m],\eta_{[m]}),u)
\]
be the nondegenerate pair which represents it. 

We will define an inverse 
\[
\lambda:\Delta^S\rightarrow \eta_*^{-1}(y)
\]
of $\pi_1\circ\omega^S$ as follows.  Let $v\in \Delta^S$.  For each $j\in 
[m]$, let $a_j\in S$ be the smallest element for which 
\[
\sum_{a\leq a_j} v_a\geq \frac{1}{k}\sum_{i=0}^j u_i
\]
Define a totally ordered set $T$ by adjoining to $S$ an immediate successor
of $a_j$, denoted $\tilde{a}_j$, for each $j$.  Define $g:T\rightarrow [m]$
by $g(b)=j$ if $\tilde{a}_{j-1}\leq b\leq a_j$. Define $f':T\rightarrow
\bar{k}$ to be $f\circ g$.  Define $h:T\rightarrow S$ by $h(b)=b$ if $b\in
S$ and $h(\tilde{a}_j)=a_j$.  Define
\[
u'_b=\left\{
\begin{array}{l}
\sum_{i=0}^j u_i-k\sum_{a<a_j}v_a \mbox{ if $b=a_j$} \\
k\sum_{a\leq a_j} v_a -\sum_{i=0}^j u_i \mbox{ if $b=\tilde{a}_j$} \\
kv_b \mbox{ otherwise}
\end{array}
\right.
\]
We define $\lambda(v)$ to be the point represented by
the pair $((f',T,h),u')$.   It is easy to check that this point is in
$\eta_*^{-1}(y)$ (this amounts to showing that
$\Upsilon((f',T,\eta_T),u')=((f,[m],\eta_{[m]}),u)$)
and that $\lambda$ is an
inverse of $\pi_1\circ\omega^S$;
this completes the proof of (2). 

Next let $e$ be a cell of $Y_k^0$ and let $\bar{e}$ be its closure.  Then 
$\eta_*^{-1}(\bar{e})$ is a finite 
union of closed cells of $Y_k^S$, and in particular it is compact.  This
together with (1) implies that
$\omega^S$ induces a homeomorphism
\[
\eta_*^{-1}(\bar{e})\rightarrow \pi_2^{-1}(\bar{e})
\]
Since the closure of each cell of $\Delta^S\times Y_k^0$ is contained in a 
set of the form $\pi_2^{-1}(\bar{e})$, it follows that $(\omega^S)^{-1}$ is 
continuous on the closure of each cell of $\Delta^S\times Y_k^0$, and from 
this it follows that $\omega^S$ is a homeomorphism.
\QED

We can now complete the proof of Theorem \ref{4.5}(a) by showing:

\begin{corollary}
\label{sat1}
$Y_k^0$ is contractible for each $k\geq 0$.
\end{corollary}

\Proof
First observe that if $A$ is a space then 
\[
\Xi_k(X_1^\b,\ldots,X_i^\b,\ldots,X_k^\b)\times A
\cong
\Xi_k(X_1^\b,\ldots,X_i^\b\times A,\ldots,X_k^\b)
\]
(because $\times A$ preserves colimits).  Thus we have
\begin{eqnarray*}
Y_k^0\times Y_j^0 &=& \Xi_k(\Delta^\b,\ldots,\Delta^\b)^0 \times Y_j^0\\
&\approx & \Xi_k(\Delta^\b,\ldots,\Delta^\b\times Y_j^0,\ldots,\Delta^\b)^0 
\\
&\approx & \Xi_k(\Delta^\b,\ldots,\Xi_j(\Delta^\b,\ldots),\ldots)^0 
\quad\mbox{by Proposition 5.3} \\
&\approx& \Xi_{k+j-1}(\Delta^\b,\ldots,\Delta^\b)^0 \quad\mbox{by Theorem
\ref{fo}} \\
&=& Y_{k+j-1}^0
\end{eqnarray*}
It therefore suffices to prove the corollary when $k=2$.  
$Y_2^0$ has the special property that the $(n-1)$-skeleton is contained in 
the closure of either of the two $n$-cells.
Corollary \ref{5.2} implies that the closure of
a cell is contractible, so
the inclusion of the $(n-1)$-skeleton is 
nullhomotopic for each $n$. Thus any map from a sphere into $Y_2^0$ is
nullhomotopic, so $Y_2^0$ is contractible.
\QED

%\newpage

\section{The structure of $\Xi_k^n(\Delta^\b,\ldots,\Delta^\b)$.}
\label{secB}

For use in Section \ref{secJ}, we prove
the analogs for $\Xi_k^n$ of the results of Section \ref{secA}.
 
Fix $n$ and 
denote the augmented cosimplicial
space $\Xi_k^n(\Delta^\b,\ldots,\Delta^\b)$ by $Z_k^\b$; thus
$\D_n(k)=\Tot(Z_k^\b)$.  Recall Notation \ref{nota}.
The elements of $Z_k^S$ are equivalence classes of pairs
$((f,T,h),u)$
with $u\in\prod_{i}\,\,\Delta^{f^{-1}(i)}$, where $f$ has complexity $\leq
n$. 

We define nondegenerate pairs exactly as in Definition \ref{nondeg}, and the
proof of Proposition \ref{5.1} goes through to show

\begin{proposition}
\label{5.1n}
Each point in $Z_k^S$ is represented by a unique nondegenerate pair. \qed
\end{proposition}

\begin{corollary}
\label{5.2n}
$Z_k^S$ is a CW complex with one cell of dimension $m+1-k$ for each
nondegenerate $(f,[m],h)$ for which $f$ has complexity $\leq n$; 
the characteristic map of the cell corresponding to $(f,[m],h)$ is a
homeomorphism from $\prod\Delta^{f^{-1}(i)}$ to the closure of the cell.
\qed
\end{corollary}

Proposition \ref{5.1n} also gives a useful relationship 
between $Z_k^\b$ and the augmented cosimplicial space $Y_k^\b$ defined in 
Section \ref{secA}:

\begin{corollary}
\label{plb}
The map $Z_k^S\rightarrow Y_k^S$ is a monomorphism for all $S$, and the 
diagram
\[
\xymatrix{
Z_k^S 
\ar[d]
\ar[r]
&
Y_k^S
\ar[d]
\\
Z_k^0
\ar[r]
&
Y_k^0
}
\]
is a pullback. \qed
\end{corollary}

Next we define 
\[
\omega:Z_k^\b\rightarrow \Delta^\b\times Z_k^0
\]
as in Section \ref{secA}: the projection of $\omega^S$ on $\Delta^S$ takes the 
equivalence class of $((f,T,h),u)$ to $v$, where 
$v_a=\frac{1}{k}\sum_{b\in h^{-1}(a)} u_b$, and the projection of $\omega^S$
on $Z_k^0$ is $\eta_*$ (see Notation \ref{eta}). The diagram
\[
\xymatrix{
Z_k^\b 
\ar[r]^-{\omega}
\ar[d]
&
\Delta^\b\times Z_k^0
\ar[d]
\\
Y_k^\b
\ar[r]^-{\omega}
&
\Delta^\b\times Y_k^0
}
\]
commutes, and this together with Corollary \ref{plb} implies

\begin{proposition}
\label{5.3n}
$\omega:Z_k^\b\rightarrow \Delta^\b\times Z_k^0$ is an isomorphism of 
augmented cosimplicial spaces. 
\end{proposition}

%\newpage

\section{Proof of Theorem \ref{6.4}(a).}

\label{secJ}

In this section we prove Theorem \ref{6.4}(a).
As motivation for the method, recall that one way to show that two 
spaces are weakly equivalent is to show that they have contractible open 
covers with the same nerve, or more generally to show that they can be 
decomposed into homotopy colimits of contractible pieces over the same 
indexing category.  We will show that the cosimplicial space 
$\Xi^n_k(\Delta^\b,\ldots,\Delta^\b)$ can be decomposed as a homotopy colimit 
of contractible cosimplicial spaces indexed over a certain category $\K^n_k$ 
considered by Berger \cite{Berger}; Berger has shown that $\C_n(k)$ is a 
homotopy colimit of contractible pieces indexed by $\K^n_k$, and from this we 
will deduce Theorem \ref{6.4}(a).

We begin by recalling some definitions from \cite{Berger} (but our notation
differs somewhat from that in \cite{Berger}).

\begin{definition}
{\rm
For each $k\ge0$, 
let $P_{2}\bar{k}$ be the set of subsets of
$\bar{k}$ that have two elements. 
}
\end{definition}

\begin{definition}  
{\rm
(a)
Let $\mathcal K_k$ be the set whose elements are pairs
$(b,T)$, where $b$ is a function from $P_{2}\bar{k}$ to the nonnegative
integers and $T$ is a total ordering of $\bar{k}$. We give 
$\mathcal K_k$ the partial order for which
$(a,S)\le (b,T)$ 
if $a(\{i,j\})\le b(\{i,j\})$ for each $\{i,j\}\in
P_{2}\bar{k}$ and $a(\{i,j\})< b(\{i,j\})$ for each $\{i,j\}$
with $i<j$ in the order $S$ but $i>j$ in the order $T$. Let
$\K^n_k$ be the subset of pairs $(b,T)$ 
such that $b\{i,j\}<n$ for each $\{i,j\}\in P_{2}\bar{k}$. The
set $\K^n_k$ inherits an order from $\K_k$.

(b) Let $\K$ denote the collection of partially ordered sets $\K_k$, $k\geq
0$, and let $\K^n$ denote the collection of partially ordered sets $\K^n_k$,
$k\geq 0$.
}
\end{definition}

It is shown in \cite{Berger} that $\K$ is an operad in the category of
partially ordered sets with the following structure maps. 
The right action of
$\Sigma_{k}$ on $\K_k$ is given by 
\[
(b,T)\rho=(b\circ\rho_{2},T\rho)
\]
where $\rho_{2}\colon P_{2}\bar{k}\to P_{2}\bar{k}$ is the function
$\rho_{2}(\{i,j\})=\{\rho(i),\rho(j)\}$ and where $i<j$ in the
total order $T\rho$ if $\rho(i)<\rho(j)$ in the total order $T$.
The operad composition
\[
\K_k\times\K_{a_{1}}\times\dots\times\K_{a_{k}}\to \K_{\Sigma a_{i}}
\] 
takes 
$((b,T);(b_1,T_{1}),\dots,(b_k,T_{k}))$ 
to the pair $(b(b_1,\dots,b_k),T(T_1,\dots,T_k))$, where
$b(b_1,\dots,b_k)$ is the function which takes $\{r,s\}$ to
\[
\left\{
\begin{array}{ll}
b_{i}(\{r,s\}) & \mbox{if $\{r,s\}\subset\ov a_{i}$} \\
b(\{i,j\}) & \mbox{if $r\in\ov a_{i}$, $s\in\ov a_{j}$ and
$i\neq j$}
\end{array}
\right.
\]
and $T(T_1,\dots,T_k)$ is the total order of 
$\coprod \ov a_{i}$ for which $r<s$ if either 
$r<s$ in the order $T_{i}$  or $r\in\ov
a_{i}$, $s\in\ov a_{j}$ and $i<j$.

Note that, for each $n$, $\K^n$ is a suboperad of $\K$.

Let us write $\N$ for the functor that takes a partially ordered 
set to the geometric realization of its nerve.  Then 
$\N\K^n$ is an operad of spaces and Berger shows (\cite[Theorem 
1.16]{Berger}) that it 
is weakly equivalent to 
the little $n$-cubes operad $\C_n$.  
To complete the proof of Theorem \ref{6.4}(a) it
therefore suffices to show that $\N\K^n$ is weakly equivalent to
$\D_n$.  We will do this by finding a homotopy colimit decomposition of the
functor-operad $\Xi^n$; first we need some definitions.

\begin{definition}
\label{nnew4}
{\rm
Let $f\colon [m]\to \bar{k}$. Then
$(b_f,T_{f})\in\K_k$ is the pair where $b_{f}(\{i,j\})$ is one
less than the complexity of the restriction of $f$ to a map
$f^{-1}(\{i,j\})\to\{i,j\}$ and where $i<j$ in the total order
$T_{f}$ if the smallest element of $f^{-1}(i)$ is less than the
smallest element of $f^{-1}(j)$.
}
\end{definition}

Recall the definition of the category $\Q_k$ (Definition \ref{Q}).

\begin{definition}
{\rm
For each pair $(b,T)\in\K_k$ let $\Q_{(b,T)}$ be the full
subcategory of $\Q_k$ whose objects are the maps $f$ with $(b_f,T_f)\le 
(b,T)$.  Let 
\[ \iota_{(b,T)}: \Q_{(b,T)}\to \Q_{k}
\]
be the inclusion functor. 
}
\end{definition}

The next definition uses the notation of Definition \ref{Kan2}.

\begin{definition}
\label{Th1}
{\rm
Let
$X^\b_{1},\dots,X^\b_{k}$ 
be augmented cosimplicial spaces.

(a) For each $(b,T)\in \K_k$,
define
$\Xi_{(b,T)}(X^\b_{1},\dots,X^\b_{k})$ to be the Kan extension
\[
\mbox{Lan}_\Phi ((X_1^\b\bar\times\cdots\bar\times X_k^\b)\circ 
\Psi\circ \iota_{(b,T)})
\]

(b) Define $\Lambda^n(X^\b_{1},\dots,X^\b_{k})$ to be 
\[
\hocolim_{\K^n_k} \Xi_{(b,T)}(X^\b_{1},\dots,X^\b_{k})
\]
}
\end{definition}

\begin{lemma}
$\Lambda^n$ is a functor-operad.
\end{lemma}

\Proof
First note that 
the natural transformations defining the functor-operad $\Xi$
restrict to natural transformations
\begin{equation}
\label{FFF2}
\sigma_*:\Xi_{(b,T)}\to \Xi_{(b,T)\sigma}\circ \sigma_\#
\end{equation}
for $\sigma\in\Sigma_{k}$ and
\begin{equation}
\label{FFF3}
\Gamma:\Xi_{(b,T)}(\Xi_{(b_{1},T_{1})},\dots,\Xi_{(b_{k},T_{k})}) \to
    \Xi_{(b(b_{1},\dots,b_{k}),T(T_{1},\dots,T_{k}))}
\end{equation}

Next recall the definition of $\hocolim$ given in \cite[Section 19.1]{H}: if 
$\A$ is a category and $F:\A\to \mbox{Top}$ is a functor then 
\[
\hocolim_\A F= F \otimes_\A U
\]
where $\otimes_\A$ denotes the coend and $U$ is the contravariant
functor $\A\to\mbox{Top}$ which takes an object $a\in\A$ to $\N(a\downarrow
\A)$.

Now let $\sigma\in\Sigma_k$ and observe that $\sigma$ induces a functor
\[
\sigma_{\downarrow}:
((b,T)\downarrow \K_k^n) \to
((b,T)\sigma\downarrow \K_k^n)
\]
We define 
\[
\sigma_*:\Lambda^n_k(X_1^\b,\ldots,X_k^\b)
\to
\Lambda^n_k(X_{\sigma(1)}^\b,\ldots,X_{\sigma(k)}^\b)
\]
to be the map induced by the collection of maps
\[
\Xi_{(b,T)}(X_1^\b,\ldots,X_k^\b) 
\times \N((b,T)\downarrow \K^n_k)
\labarrow{\sigma_*\times\N(\sigma_{\downarrow})}
\Xi_{(b,T)\sigma}(X_{\sigma(1)}^\b,\ldots,X_{\sigma(k)}^\b)
\times
\N((b,T)\sigma\downarrow \K^n_k)
\]

Finally, we define the structural map
\[
\Gamma:\Lambda^n_k(\Lambda^n_{j_1},\ldots,\Lambda^n_{j_k})
\to
\Lambda^n_{j_1+\cdots+j_k}
\]
to be the map induced by the collection of maps
\begin{multline}
\notag
\Xi_{(b,T)}(\Xi_{(b_{1},T_{1})},\dots,\Xi_{(b_{k},T_{k})}) \times
\N((b,T)\downarrow \K^n_k)\times
\prod_{i=1}^k \N((b_i,T_i)\downarrow \K^n_{j_i})
\labarrow{\cong}
\\
\Xi_{(b,T)}(\Xi_{(b_{1},T_{1})},\dots,\Xi_{(b_{k},T_{k})}) \times
\N\Bigl(((b,T)\downarrow \K^n_k) \times 
\prod_{i=1}^k ((b_i,T_i)\downarrow \K^n_{j_i})\Bigr)
\labarrow{\Gamma \times \N(\gamma_{\downarrow})}
\\
\Xi_{(b(b_{1},\dots,b_{k}),T(T_{1},\dots,T_{k}))}
\times
\N((b(b_{1},\dots,b_{k}),T(T_{1},\dots,T_{k}))\downarrow\K^n_{j_1+\cdots+j_k})
\end{multline}
where $\gamma_{\downarrow}$ is induced by the composition map
\[
\gamma:\K^n_k\times\prod_{i=1}^k \K^n_{j_i}
\to
\K^n_{j_1+\cdots+j_k}
\]
of the Cat-operad $\K^n$.
\QED

Now let $\B_n$ be the operad obtained by
applying Proposition \ref{3.2} with $\F=\Lambda^n$ and $A=\Delta^\b$.
To complete the proof of Theorem \ref{6.4}(a) it remains to show:

\begin{lemma}
\label{WW1}
(a) There is a weak equivalence of operads
\[
\B_n\to \D_n
\]

(b) There is a weak equivalence of operads 
\[
\B_n \to \N\K^n
\]
\end{lemma}

For the proof of part (a), we first observe that $\Q^n_k$ is the union 
of $\Q_{(b,T)}$ for $(b,T)\in \K^n_k$; it follows that
\[
\Xi^n_k(X^\b_{1},\dots,X^\b_{k})=
\colim_{\K^n_k} \Xi_{(b,T)}(X^\b_{1},\dots,X^\b_{k})
\]
for all $X^\b_1,\ldots,X^\b_k$.
The projection from hocolim to colim gives a map of functor-operads
\[
\Lambda^n\to \Xi^n
\]
and an induced map of the associated operads:
\[
\phi:\B_n\to \D_n
\]
We need to show that for each $k\geq 0$ the map
\[
\phi(k):\B_n(k)\to \D_n(k)
\]
is a weak equivalence of spaces.  Recall from Proposition \ref{5.3n}
that there is an isomorphism of cosimplicial spaces
\[
\Xi^n_k(\Delta^\b,\ldots,\Delta^\b)
\cong
\Delta^\b \times
\Xi^n_k(\Delta^\b,\ldots,\Delta^\b)^0
\]
The proof of Proposition \ref{5.3n} shows that for each $(b,T)$ we have
an isomorphism of cosimplicial spaces
\[
\Xi_{(b,T)}(\Delta^\b,\dots,\Delta^\b)
\cong \Delta^\b \times
\Xi_{(b,T)}(\Delta^\b,\dots,\Delta^\b)^0
\]
It follows that we have homeomorphisms
\[
\D_n(k) \approx 
\Tot(\Delta^\b) \times
\colim_{\K^n_k} 
\Xi_{(b,T)}(\Delta^\b,\ldots,\Delta^\b)^0
\]
and
\begin{equation}
\label{FFF1}
\B_n(k) \approx 
\Tot(\Delta^\b) \times
\hocolim_{\K^n_k} 
\Xi_{(b,T)}(\Delta^\b,\ldots,\Delta^\b)^0
\end{equation}
Thus it suffices to show that the map
\[
\hocolim_{\K^n_k}
\Xi_{(b,T)}(\Delta^\b,\ldots,\Delta^\b)^0 \to 
\colim_{\K^n_k}
\Xi_{(b,T)}(\Delta^\b,\ldots,\Delta^\b)^0
\]
is a weak equivalence, and this follows from a standard fact about homotopy
colimits \cite[Theorem 19.9.1]{H}; the ``Reedy cofibrancy'' condition needed
for \cite[Theorem 19.9.1]{H} is satisfied in our case because the map
\[
\bigcup_{(b',T')<(b,T)} \Xi_{(b',T')}(\Delta^\b,\dots,\Delta^\b)^0
\to \Xi_{(b,T)}(\Delta^\b,\dots,\Delta^\b)^0
\]
is the inclusion of a sub-CW-complex (cf.\ Corollary \ref{5.2n}).

Next we prove part (b).  Consider the map
\[
\psi(k):
\B_n(k)  = \Hom(\Delta^\b,\Lambda^n_k(\Delta^\b,\ldots,\Delta^\b))
\to
\Hom(\Delta^\b,\Lambda^n_k(*,\ldots,*))
=
\N\K^n_k
\]
where the arrow is induced by the projection $\Delta^\b\to *$ and the second
equality follows from the fact that $\Lambda^n_k(*,\ldots,*)$ is the constant 
cosimplicial space with value $\N\K^n_k$.  It is easy to check 
that the collection $\{\psi(k)\}$ is an operad map; it remains to show that 
each $\psi(k)$ is a weak equivalence.   
Using equation \eqref{FFF1} and the fact that
$\Tot(\Delta^\b)$ is contractible 
it suffices to show that the map
\[
\hocolim_{\K^n_k} \Xi_{(b,T)}(\Delta^\b,\dots,\Delta^\b)^0
\to
\hocolim_{\K^n_k} \Xi_{(b,T)}(*,\ldots,*)^0
=\hocolim_{\K^n_k} *
\]
(where the arrow is induced by $\Delta^\b\to *$) is a weak equivalence; and
this is a consequence of \cite[Remark 18.5.4]{H} and the following lemma.

\begin{lemma}
\label{nnnew}
For each choice of $b$ and $T$ 
the space $\Xi_{(b,T)}(\Delta^\b,\dots,\Delta^\b)^0$ is
weakly equivalent to a point.
\end{lemma}

\Proof
The proof is by induction on $k$.  

Since the map \eqref{FFF2}
is an isomorphism we may assume that $T$ is the standard total order on
$\bar{k}$.

For each $f:[m]\to \bar{k}$ we define $c(f):[m+1]\to \bar{k}$ to be the 
function which 
takes 1 to 1 and $p$ to $f(p-1)$ if $p>1$.
This construction gives a functor, also called $c$, from $\Q_{(b,T)}$ 
to itself.

Next let $C:\Delta\to\Delta$ be the functor which takes $[m]$ to $[m+1]$ and
takes a morphism $h:[m]\to[n]$ to the morphism $C(h):[m+1]\to [n+1]$ defined
by
\[
C(h)(p)=\left\{
\begin{array}{ll}
0 & \mbox{if $p=0$} \\
h(p-1)+1 & \mbox{if $p>0$}
\end{array}
\right.
\]
We can define a map
\[
\alpha:
\Xi_{(b,T)}(\Delta^\b,\ldots,\Delta^\b)^0
\to
\Xi_{(b,T)}(\Delta^\b\circ C,\Delta^\b,\ldots,\Delta^\b)^0
\]
as follows: if $f:[p]\to\bar{k}$ and $u_i\in \Delta^{f^{-1}(i)}$ for $1\leq
i\leq k$, let $\alpha$ take the equivalence class of 
$(f,u_1,\ldots,u_n)$ to that of $(f,d^0 u_1, u_2,\ldots,u_k)$.
We can also define a map 
\[
\beta:
\Xi_{(b,T)}(\Delta^\b\circ C,\Delta^\b,\ldots,\Delta^\b)^0
\to
\Xi_{(b,T)}(\Delta^\b,\ldots,\Delta^\b)^0
\]
by letting $\beta$ take the equivalence class of 
$(f,u_1,\ldots,u_n)$ to that of 
$(cf,u_1,\ldots,u_k)$.  It is easy to check that $\alpha$ and $\beta$ are
well-defined and that $\beta\circ\alpha$ is the identity; that is, 
$\Xi_{(b,T)}(\Delta^\b,\ldots,\Delta^\b)^0$
is a retract of 
$\Xi_{(b,T)}(\Delta^\b\circ C,\Delta^\b,\ldots,\Delta^\b)^0$. It therefore
suffices to show that the latter is weakly equivalent to a point.

But $\Delta^\b\circ C$ is isomorphic to the degreewise cone on $\Delta^\b$, and
in particular there is a homotopy equivalence of cosimplicial spaces from 
$\Delta^\b\circ C$ to a point, so we have
\[
\Xi_{(b,T)}(\Delta^\b\circ C,\Delta^\b,\ldots,\Delta^\b)^0
\simeq
\Xi_{(b,T)}(*,\Delta^\b,\ldots,\Delta^\b)^0
\]
and an inspection of Definition \ref{Th1} shows that 
\[
\Xi_{(b,T)}(*,\Delta^\b,\ldots,\Delta^\b)
\cong
\Xi_{(b',T')}(\Delta^\b,\ldots,\Delta^\b)
\]
where $b'$ is the restriction of $b$ to $P_2(\bar{k}-\{1\})$ and $T'$ is the
restriction of $T$ to $\bar{k}-\{1\}$.  The inductive hypothesis shows that 
$\Xi_{(b',T')}(\Delta^\b,\ldots,\Delta^\b)^0$ is weakly equivalent to a point,
and this concludes the proof.
\QED

%\newpage

\section{A homotopy-invariant version of Tot.}

\label{sec8}

First recall (\cite[Theorem 11.6.1]{H}) that there is a model category 
structure for cosimplicial spaces in which the weak equivalences are the 
degreewise weak equivalences and the fibrations are the degreewise fibrations.
In particular, every object is fibrant.  Let $\tilde{\Delta}^\b$ be any
cofibrant resolution of $\Delta^\b$ with respect to this model structure.

\begin{definition}
{\rm
Let $X^\b$ be a cosimplicial space.  $\TOT(X^\b)$ is defined to be 
$\Hom(\tilde{\Delta}^\b,X^\b)$.
}
\end{definition}

Since every cosimplicial space is fibrant, a weak equivalence $X^\b\to Y^\b$
always induces a weak equivalence $\TOT(X^\b)\to\TOT(Y^\b)$.

\begin{definition}
{\rm
Let $\tilde{\D}_n$ be the operad obtained by applying Proposition \ref{3.2} 
with $\F=\Xi^n$ and $A=\tilde{\Delta}^\b$ 
}
\end{definition}

Our goal in this section is to prove the analog of Theorem \ref{6.4}.

\begin{theorem}
\label{nnew1}
{\rm (a)}
$\tilde{\D}_n$ is weakly equivalent in the category of operads to
$\C_n$.

{\rm (b)}
If $X^\b$ is an algebra over $\Xi^n$ then $\tilde{\D}_n$ acts on
$\TOT(X^\b)$.
\end{theorem}

\begin{remark}
{\rm
Note that this theorem includes the analog of Theorem \ref{4.5} as a special
case.  Moreover, the proof we will give can easily be modified to prove the
analog of Proposition \ref{2.3}.
}
\end{remark}

Before beginning the proof we give some background information which is of
interest in its own right.  

We begin with a more explicit description of 
$\Xi_k^n(\tilde{\Delta}^\b,\ldots,\tilde{\Delta}^\b)$.
First observe that $\tilde{\Delta}^\b$ is a Reedy-cofibrant cosimplicial space
by \cite[Proposition 15.6.3(2)]{H}.  It follows that
\[
\tilde{\Delta}^\b\times\dots\times\tilde{\Delta}^\b
\]
is a Reedy-cofibrant multicosimplicial space.
Next observe that we can extend the definition of $\Xi_k^n$ to 
$k$-fold multicosimplicial spaces $Y^{\b,\ldots,\b}$ 
by replacing 
$X_1^\b\bar\times\cdots\bar\times X_k^\b$ in 
Definition \ref{Kan3} by $Y^{\b,\ldots,\b}$.  
Now fix a finite totally ordered set $S$.
For each $m\geq k-1$ let $I_m$ be the set of nondegenerate diagrams
$(f,[m],h)$ where $f$ has complexity $\leq n$
(see Notation \ref{nota} and Definition \ref{nondeg}(a); 
note that there cannot be any nondegenerate $(f,[m],h)$ if $m< k-1$).
Recall the definition of latching object (\cite[Definition 15.2.5]{H}).

\begin{lemma}
\label{MLK1}
Let $Y^{\b,\ldots,\b}$ be a Reedy-cofibrant $k$-fold multicosimplicial space.
Define 
\[
V_{k-1}=\coprod_{(f,[m],h)\in I_{k-1}} Y^{0,\ldots,0}
\]
and define $V_m$ inductively for $m\geq k$ by the pushout diagram
\[
\xymatrix{
\coprod L_f \
\ar[d]
\ar[r]
&
\ V_{m-1}
\ar[d]
\\
\coprod Y^{f^{-1}(1),\ldots,f^{-1}(k)} \
\ar[r]
&
V_m
}
\]
where the coproducts are taken over $(f,[m],h)\in I_m$
and $L_f$ is the latching object 
\[
L_{f^{-1}(1),\ldots,f^{-1}(k)}(Y^{\b,\ldots,\b}).
\]
Then each $V_{m-1}\to V_m$ is a cofibration, and
\[
\Xi_k^n(Y^{\b,\ldots,\b})^S=\bigcup_{m\geq k-1} V_m.
\]
\end{lemma}

\Proof
This follows by the proof of 
Corollary \ref{5.2n}. 
\QED

Next we give a homotopy-invariance property for $\Xi^n_k$.

\begin{lemma}
\label{MLK2}
Let $Y^{\b,\ldots,\b}\to Z^{\b,\ldots,\b}$ be a weak equivalence of 
Reedy-cofibrant $k$-fold multicosimplicial spaces.  Then the induced map
\[
\Xi^n_k(Y^{\b,\ldots,\b})\to
\Xi^n_k(Z^{\b,\ldots,\b})
\]
is a weak equivalence.
\end{lemma}

\Proof
This is an easy consequence of Lemma \ref{MLK1}.
\QED

\begin{remark}
\label{MLK3}
{\rm
The analog of Lemma \ref{MLK2} for the functors $\Xi_{(b,T)}$ defined in
Section \ref{secJ} is also true, with the same proof.
}
\end{remark}

Now we turn to the proof of Theorem \ref{nnew1}.
Part (b) is immediate from Proposition \ref{3.4}.  For part (a),
recall the functor-operad $\Lambda^n$ and the operad $\N\K^n$
defined in Section \ref{secJ}.  Let $\tilde{\B}_n$ be the operad obtained by
applying Proposition \ref{3.2} with $\F=\Lambda^n$ and $A=\tilde{\Delta}^\b$.
It suffices to show

\begin{lemma}
\label{nnew2}
(i) There is a weak equivalence of operads
\[
\tilde{\B}_n\to \tilde{\D_n}
\]

(ii) There is a weak equivalence of operads
\[
\tilde{\B}_n \to \N\K^n
\]
\end{lemma}

\noindent{\bf Proof of Lemma \ref{nnew2}.\ }
For part (i), note that
the projection from hocolim to colim gives a map of functor-operads
\[
\Lambda^n\to \Xi^n
\]
and an induced map of the associated operads:
\[
\tilde{\phi}:\tilde{\B}_n\to \tilde{\D}_n
\]
We need to show that for each $k\geq 0$ the map
\[
\tilde{\phi}(k):\tilde{\B}_n(k)\to \tilde{\D}_n(k)
\]
is a weak equivalence of spaces.  Because $\TOT$ is homotopy-invariant, it
suffices to show that the map
\[
\Lambda^n_k(\tilde{\Delta}^\b,\ldots,\tilde{\Delta}^\b)^S \to
\Xi^n_k(\tilde{\Delta}^\b,\ldots,\tilde{\Delta}^\b)^S
\]
is a weak equivalence for each $S$.
Consider the commutative diagram
\[
\xymatrix{
\Lambda^n_k(\tilde{\Delta}^\b,\ldots,\tilde{\Delta}^\b)^S
\ar[r]
\ar[d]
&
\Xi^n_k(\tilde{\Delta}^\b,\ldots,\tilde{\Delta}^\b)^S
\ar[d]
\\
\Lambda^n_k({\Delta}^\b,\ldots,{\Delta}^\b)^S
\ar[r]
&
\Xi^n_k({\Delta}^\b,\ldots,{\Delta}^\b)^S
}
\]
where the vertical maps are induced by the projection
$\tilde{\Delta}^\b\to\Delta^\b$.
We have shown in the proof of Lemma \ref{WW1}(a) that the lower horizontal map
is a weak equivalence.  
The second vertical map is a weak equivalence by Lemma \ref{MLK2}, and the
first vertical map is a weak equivalence by Remark \ref{MLK3} and passage to
hocolim.
This completes the proof of part (i).

For part (ii),
let $D$ be the 0-th space of $\tilde{\Delta}^\b$.  
If 
$C^\b$ is any constant cosimplicial space there is a canonical homeomorphism
\[
\TOT(C^\b)\approx \Map(D,C^0)
\]
Now $\Lambda^n_k(*,\ldots,*)$ is the constant cosimplicial space
with value $\N\K^n_k$, so we have 
\[
\TOT(\Lambda^n_k(*,\ldots,*)) \approx \Map(D, \N\K^n_k).
\]
Thus $\TOT(\Lambda^n_k(*,\ldots,*))$ is an operad weakly 
equivalent to $\N\K^n$.  It therefore suffices to show that the map
\[
\tilde{\B}_n(k)  = 
\TOT(\Lambda^n_k(\tilde{\Delta}^\b,\ldots,\tilde{\Delta}^\b))
\to
\TOT(\Lambda^n_k(*,\ldots,*))
\]
induced by the projection $\tilde{\Delta}^\b\to *$ is a weak equivalence for 
each $k$.   Since $\TOT$ preserves weak equivalences, it suffices to show that
the map
\[
\hocolim_{\K^n_k} \Xi_{(b,T)}(\tilde{\Delta}^\b,\ldots,\tilde{\Delta}^\b)^S
\to
\hocolim_{\K^n_k} \Xi_{(b,T)}(*,\ldots,*)^S =
\hocolim_{\K^n_k} *
\]
is a weak equivalence for every $k$ and $S$, and this in turn follows from 
Lemma \ref{nnnew} and Remark \ref{MLK3}.

\end{document}